\documentclass[11pt,a4paper]{article}
\usepackage[latin1]{inputenc}
\usepackage{amsmath}
\usepackage{amsfonts}
\usepackage{amssymb}
\usepackage{amsthm}
\usepackage{stmaryrd}

\DeclareMathOperator{\E}{\mathbb{E}}
\DeclareMathOperator{\N}{\mathbb{N}}
\DeclareMathOperator{\R}{\mathbb{R}}
\DeclareMathOperator{\C}{\mathbb{C}}
\DeclareMathOperator{\PR}{\mathbb{P}}
\DeclareMathOperator{\BB}{\mathcal{B}}
\DeclareMathOperator{\CC}{\mathcal{C}}
\DeclareMathOperator{\FF}{\mathcal{F}}
\DeclareMathOperator{\HH}{\mathcal{H}}
\DeclareMathOperator{\MM}{\mathcal{M}}
\DeclareMathOperator{\NN}{\mathcal{N}}
\DeclareMathOperator{\OO}{\mathcal{O}}
\DeclareMathOperator{\PP}{\mathcal{P}}
\DeclareMathOperator{\Sz}{\mathcal{S}}
\DeclareMathOperator{\1}{\mathbf{1}}
\DeclareMathOperator{\Reel}{\mathrm{Re}}
\DeclareMathOperator{\Imag}{\mathrm{Im}}
\DeclareMathOperator{\Int}{\mathrm{Int}}
\DeclareMathOperator{\Adh}{\mathrm{Clo}}
\DeclareMathOperator{\Var}{\mathrm{Var}}
\DeclareMathOperator{\Id}{\mathrm{Id}}
\DeclareMathOperator{\Tr}{\mathrm{Tr}}
\DeclareMathOperator{\rg}{\mathrm{rank}}
\DeclareMathOperator{\diag}{\mathrm{diag}}
\DeclareMathOperator{\sgn}{\mathrm{sign}}
\DeclareMathOperator{\musc}{\mu_{\mathrm{sc}}}
\newcommand{\mump}{\mu_{\mathrm{MP},c}}
\newcommand{\dKS}{d_{\mathrm{KS}}}
\newcommand{\gbar}{\underline{g}}
\newcommand{\Xhat}{\widehat{X}}
\newcommand{\Yhat}{\widehat{Y}}

\newtheorem{Proposition}{Proposition}[section]
\newtheorem{Theorem}[Proposition]{Theorem}
\newtheorem*{Theorem*}{Theorem}
\newtheorem{Lemma}[Proposition]{Lemma}

\theoremstyle{definition}
\newtheorem{Definition}[Proposition]{Definition}
\newtheorem*{Remark}{Remark}
\newtheorem*{Remarks}{Remarks}

\newtheorem*{Proof}{Proof}

\title{Asymptotic Freeness for Rectangular Random Matrices and Large Deviations for Sample Covariance Matrices With Sub-Gaussian Tails}
\author{Benjamin \textsc{Groux} \footnote{Universit\'e de Versailles Saint-Quentin-en-Yvelines, Laboratoire de Math\'ematiques de Versailles, 45 avenue des \'Etats-Unis, F-78035 \textsc{Versailles Cedex}. E-mail: \texttt{benjamin.groux@uvsq.fr} } }
\date{\today}

\begin{document}

\maketitle

\begin{abstract}
We establish a large deviation principle for the empirical spectral measure of a sample covariance matrix with sub-Gaussian entries, which extends Bordenave and Caputo's result for Wigner matrices having the same type of entries \cite{BC}. To this aim, we need to establish an asymptotic freeness result for rectangular free convolution, more precisely, we give a bound in the subordination formula for information-plus-noise matrices.\\
\end{abstract}

\textbf{AMS 2010 Classification Subject.} 60B20, 46L54, 60F10, 15A18.\\

\textbf{Key words.} Random matrices; Large deviations; Free convolution; Subordination property; Spectral measure; Stieltjes transform; Information-plus-noise model.

\tableofcontents

\section{Introduction}

Throughout this paper, $\PP(E)$ will denote the set of probability measures on a space $E$, $\MM_{n,p}(\R)$ (resp. $\MM_{n,p}(\C)$) the set of $n \times p$ real (resp. complex) matrices, $\HH_n(\C)$ the set of $n \times n$ Hermitian matrices, $A^t$ (resp. $A^*$) the transpose (resp. transconjugate) of a matrix $A$, and $\Tr(A)$ its trace. Besides, for a random variable $X$, $\mathring{X}$ denotes the centred variable $X-\E(X)$. Finally, for two real numbers $x,y$, we denote by $x \wedge y$ the minimum of $x$ and $y$.

\subsection{Large deviation results in random matrix theory}

Let us first recall some basic facts in random matrix theory (RMT). A key object in RMT is the \emph{empirical spectral measure} of a matrix $A \in \HH_n(\C)$, namely the probability measure on $\R$ defined by
\begin{displaymath}
\mu_A = \frac1n \sum_{k=1}^n \delta_{\lambda_k(A)} \, ,
\end{displaymath}
where $\lambda_1(A),\ldots,\lambda_n(A)$ denote the eigenvalues of $A$.

It is well known (cf. \cite{W}) that if $X$ is a \emph{Wigner matrix}, i.e. $X \in \HH_n(\C)$ and the families of centred independent and identically distributed (i.i.d.) random variables $(X_{j,j})_{1 \le j \le n}$, $(X_{j,k})_{1 \le j < k \le n}$ are independent, and if the variance $\Var(X_{1,2}) = \E |X_{1,2} - \E(X_{1,2})|^2$ equals 1, then almost surely, the spectral measure $\mu_{X/\sqrt{n}}$ converges weakly towards the \emph{semicircular distribution} $\musc$, i.e. for any bounded continuous $f : \R \to \R$,
\begin{displaymath}
\lim_{n \to +\infty} \int_{\R} f \, d\mu_{X/\sqrt{n}} = \int_{\R} f \, d\musc \, .
\end{displaymath}
The semicircular distribution $\musc$ is the probability measure on $\R$ defined by
\begin{displaymath}
d\musc(x) = \frac{1}{2\pi} \sqrt{4-x^2} \1_{[-2,2]}(x) \, dx \, .
\end{displaymath}

In the case of a \emph{sample covariance matrix}, i.e. a matrix $XX^*$ with $X \in \MM_{n,p}(\C)$ having centred i.i.d. entries, if $\Var(X_{1,1})=1$, then almost surely, the spectral measure $\mu_{XX^*/p}$ converges weakly towards the \emph{Marcenko-Pastur distribution} $\mump$ with ratio $c$ as $n,p \to +\infty$ with $\frac{n}{p} \to c \in (0,+\infty)$ (cf. \cite{MP}). This probability measure on $\R$ is defined by
\begin{displaymath}
d\mump(x) = \max \left( 1-\frac1c , 0 \right) \delta_0 + \frac{\sqrt{(b_c-x)(x-a_c)}}{2\pi xc} \1_{[a_c,b_c]}(x) \, dx
\end{displaymath}
with $a_c=(1-\sqrt{c})^2$ and $b_c=(1+\sqrt{c})^2$.

For these two models in which the empirical spectral measure converges, we can investigate the speed of convergence and more particularly large deviation principles.\\

We recall from \cite{DZ} that a sequence of random variables $(Z_n)_{n \ge 1}$ with values in a topological space $(E,\OO)$ with $\sigma$-Borel field $\BB$ satisfies the \emph{large deviation principle} (LDP) with speed $v$ and rate function $I$ in the topology $\OO$ if
\begin{itemize}
\item $I : E \to [0,+\infty]$ is a lower semi-continuous function, i.e. the level set $\{x \in E \ | \ I(x) \le t\}$ is closed for every $t \ge 0$,
\item $v : \N \to (0,+\infty)$ admits a limit equal to $+\infty$,
\item for all $B \in \BB$,
\begin{multline*}
-\inf_{x \in \Int(B)} I(x) \le \liminf_{n \to +\infty} \frac{1}{v(n)} \log \PR(Z_n \in B) \\
\le \limsup_{n \to +\infty} \frac{1}{v(n)} \log \PR(Z_n \in B) \le -\inf_{x \in \Adh(B)} I(x)
\end{multline*}
where $\Int(B)$ and $\Adh(B)$ denote resp. the interior and the closure of $B$.
\end{itemize}
We also recall that the rate function $I$ is said to be \emph{good} if the level set $\{x \in E \ | \ I(x) \le t\}$ is compact for every $t \ge 0$.\\

In \cite{BAG}, Ben Arous and Guionnet proved that if $X$ is in the GUE, i.e. $X$ is a Wigner matrix and $X_{1,1}$ (resp. $X_{1,2}$) has law $\NN(0,1)$ (resp. $\NN_2 \left( 0,\frac12 I_2 \right)$), then $\mu_{X/\sqrt{n}}$ satisfies a LDP in $\PP(\R)$ at speed $n^2$ with the rate function
\begin{displaymath}
I(\mu) = \frac12 \int x^2 \, d\mu(x) - \iint \log|x-y| \, d\mu(x) d\mu(y) - \frac34 \, .
\end{displaymath}
This result was extended to LUE matrices, i.e. sample covariance matrices $XX^*$ where $X$ has standard Gaussian entries, by Hiai and Petz (see \cite{HP}). Note that in fact, these two LDPs do not concern only Gaussian matrices but also more general unitarily invariant models. They strongly rely on the fact that for the considered models, the joint distribution of the eigenvalues has an explicit form, which is also the case in \cite{ESS}.

In \cite{BC}, Bordenave and Caputo managed to obtain a LDP for Wigner matrices in another case, where the distribution of the $X_{j,k}$'s has sub-Gaussian tails. This is remarkable because here the joint distribution of the eigenvalues is unknown. Let us recall their result.

\begin{Definition} \label{S_alpha(a)}
For $\alpha>0$ and $a \in (0,+\infty]$, we denote by $\Sz_{\alpha}(a)$ the class of complex random variables $Z$ such that
\begin{equation} \label{hypo1}
\lim_{t \to +\infty} -t^{-\alpha} \log \PR(|Z| \ge t) = a
\end{equation}
and such that $|Z|$ and $Z/|Z|$ are independent for large values of $|Z|$, i.e. there exist $t_0>0$ and a probability measure $\vartheta_a$ on the unit circle $\mathbb{S}^1$ such that for all $t \ge t_0$ and all measurable sets $U \subset \mathbb{S}^1$, we have
\begin{displaymath}
\PR(Z/|Z| \in U \ \cap \ |Z| \ge t) = \vartheta_a(U) \PR(|Z| \ge t) \, .
\end{displaymath}
In particular, a real random variable $Z$ belongs to $\Sz_{\alpha}(a)$ if it satisfies (\ref{hypo1}) and there exist $t_0>0$ and a probability measure $\vartheta_a$ on $\{-1,1\}$ such that for all $t \ge t_0$ and all $U \subset \{-1,1\}$, we have
\begin{equation} \label{hypo2bis}
\PR(|Z| \ge t \ \cap \ \sgn(Z) \in U) = \vartheta_a(U) \PR(|Z| \ge t) \, .
\end{equation}
\end{Definition}

Note that the first hypothesis implies that a random variable in $\Sz_{\alpha}(a)$ has finite moments of all orders.
\begin{Theorem}[see {\cite[Theorem 1.1]{BC}}] \label{LDP_BC}
Let $X$ be a Wigner matrix with $X_{1,2} \in \Sz_{\alpha}(a)$ and $X_{1,1} \in \Sz_{\alpha}(b)$ for some $\alpha \in (0,2)$ and $a,b \in (0,+\infty]$. Then the spectral measure $\mu_{X/\sqrt{n}}$ satisfies the LDP with speed $n^{1+\alpha/2}$ and good rate function
\begin{displaymath}
J(\mu) = \left\{ \begin{array}{ll}
\Phi(\nu) & \textrm{ if there exists } \nu \in \PP(\R) \textrm{ such that } \mu = \musc \boxplus \nu \\
+\infty & \textrm{ otherwise}
\end{array} \right.
\end{displaymath}
where $\Phi : \PP(\R) \to [0,+\infty]$ is a good rate function (see \cite{BC} for further details) and $\boxplus$ denotes the free convolution (see Section \ref{deformed_models}).
\end{Theorem}
Let us make a few remarks about this result. Roughly speaking, after random matrix considerations, the proof of Theorem \ref{LDP_BC} consists in proving a LDP for some random graphs associated to the Wigner matrix $X$. Therefore, the rate function $\Phi$ expresses as the supremum of functions of probability measures on graphs and it can not be computed in general. However, in some particular cases, it is possible to compute $\Phi(\nu)$. For example, if $\nu$ is a symmetric distribution on $\R$, $b<\infty$ and the support of $\vartheta_b$ is $\{-1,1\}$, then we have
\begin{displaymath}
\Phi(\nu) = \left( \frac{a}{2} \wedge b \right) m_{\alpha}(\nu) \, ,
\end{displaymath}
where $m_{\alpha}(\nu)$ denotes the $\alpha$-th moment of $\nu$.

Theorem \ref{theorem2} below will extend Theorem \ref{LDP_BC} to sample covariance matrices $XX^*$ with $X_{1,1} \in \Sz_{\alpha}(a)$ for some $\alpha \in (0,2)$, $a \in (0,+\infty]$. Note that to simplify, we will assume that $X$ is a \underline{real} random matrix.

Let us mention here that LDPs for the top eigenvalue of Wigner matrices have also been obtained in Ben Arous and Guionnet's setting, see \cite[p. 81]{AGZ}, and for the model introduced by Bordenave and Caputo in \cite{A}.

\subsection{Deformed matrix models} \label{deformed_models}

After understanding the behaviour of the spectral measure of Wigner matrices or sample covariance matrices, the question of deformations of these models has been investigated. Several types of deformations have been studied, the main ones being matrices of the type $X+A$ with $A \in \HH_n(\C)$ (additive deformation), $\Sigma^{1/2}XX^*\Sigma^{1/2}$ with $\Sigma \in \HH_n(\C)$ definite positive (multiplicative deformation) or $(X+A)(X+A)^*$ with $A \in \MM_{n,p}(\C)$ (information-plus-noise model).

A tool to study the spectral measure of a deformation is free probability, and more particularly free convolutions. Let us recall their definitions.

\begin{Theorem}[see \cite{V}]
Let $A$, $B$ be two independent $n \times n$ Hermitian random matrices such that
\begin{itemize}
\item either $A$ or $B$ is unitarily invariant, i.e. for $M = A$ or $B$, for any unitary $U \in \MM_n(\C)$, $UMU^*$ has the same law as $M$,
\item $\mu_A$ and $\mu_B$ converge weakly in probability to some distributions $\mu_1$ and $\mu_2$ on $\R$ as $n \to +\infty$.
\end{itemize}
Then, as $n \to +\infty$, the spectral measure $\mu_{A+B}$ converges weakly in probability to a deterministic distribution depending only on $\mu_1$ and $\mu_2$. This distribution is called the \emph{free (additive) convolution} of $\mu_1$ and $\mu_2$, and is denoted by $\mu_1 \boxplus \mu_2$.
\end{Theorem}

A similar result also exists for the singular values of the sum of two rectangular matrices and it is due to Benaych-Georges. The \emph{empirical singular value distribution} of a matrix $A \in \MM_{n,p}(\C)$ is the probability measure on $\R_+$ defined by
\begin{displaymath}
\nu_A = \frac{1}{n \wedge p} \sum_{k=1}^{n \wedge p} \delta_{\sigma_k(A)} \, ,
\end{displaymath}
where $\sigma_1(A),\ldots,\sigma_{n \wedge p}(A)$ denote the singular values of $A$, i.e. the square roots of the eigenvalues of the positive matrix $AA^*$ (resp. $A^*A$) if $n \le p$ (resp. $n \ge p$).

\begin{Theorem}[see {\cite[Theorem 3.13]{BG}}]
Let $A$, $B$ be two independent $n \times p$ random matrices such that
\begin{itemize}
\item either $A$ or $B$ is bi-unitarily invariant, i.e. for $M = A$ or $B$, for any unitary matrices $U \in \MM_n(\C)$ and $V \in \MM_p(\C)$, $UMV$ has the same law as $M$,
\item $\nu_A$ and $\nu_B$ converge weakly in probability to some distributions $\mu_1$ and $\mu_2$ on $\R_+$ as $n,p \to +\infty$ with $\frac{n}{p} \to c \in (0,+\infty)$.
\end{itemize}
Then, as $n \to +\infty$, the singular value distribution $\nu_{A+B}$ converges weakly in probability to a deterministic distribution depending only on $\mu_1$, $\mu_2$ and $c$. This distribution is called the \emph{rectangular free convolution} with ratio $c$ of $\mu_1$ and $\mu_2$, and is denoted by $\mu_1 \boxplus_c \mu_2$.
\end{Theorem}

Free convolutions can be characterized in terms of another key object in RMT, Stieltjes transform. For a probability measure $\mu$ on $\R$, we call the \emph{Stieltjes transform} of $\mu$ the function $G_{\mu} : \C \setminus \R \to \C$ defined by
\begin{displaymath}
G_{\mu}(z) = \int_{\R} \frac{1}{z-x} \, d\mu(x)
\end{displaymath}
for all $z \in \C \setminus \R$. The following properties are obvious:
\begin{displaymath}
|G_{\mu}(z)| \le \frac{1}{|\Imag z|}
\end{displaymath}
and
\begin{displaymath}
|G_{\mu}(z) - G_{\mu}(z')| \le \frac{|z-z'|}{|\Imag z|.|\Imag z'|} \, .
\end{displaymath}
We will use them implicitly in this paper.

Note that the notion of Stieltjes transform is related to the resolvent one, since for a matrix $A \in \HH_n(\C)$, we have $G_{\mu_A}(z) = \frac1n \Tr((zI_n-A)^{-1})$. Useful properties of resolvents we will use in this paper are gathered in Appendix \ref{resolvents}.\\

Stieltjes transform allows to express subordination relations for free convolutions. To state these relations, we need some additional notations. For $\mu \in \PP(\R)$, we denote by $\mu^2$ the distribution of $X^2$ when $X$ has law $\mu$. Similarly, for $\mu \in \PP(\R_+)$, we denote by $\sqrt{\mu}$ the symmetrization of the distribution $\nu$ of $\sqrt{X}$ when $X$ has law $\mu$, i.e. the symmetric distribution on $\R$ defined by $\sqrt{\mu}(B) = \frac{\nu(B)+\nu(-B)}{2}$ for all borelians $B$. We have the following subordination formulas, the first is due to Biane (cf. \cite{Biane}) and the second is obtained from Dozier and Silverstein's work \cite{DS} and a paper by Benaych-Georges (cf. \cite{BG}).

\begin{Proposition}
\begin{itemize}
\item Let $\mu \in \PP(\R)$ and $\nu = \mu \boxplus \mu_{sc}$. We have
\begin{equation} \label{subordination_square}
G_{\nu}(z) = G_{\mu} \left( z-G_{\nu}(z) \right) \, .
\end{equation}
\item Let $\mu \in \PP(\R_+)$, $c>0$ and $\nu = \left( \sqrt{\mu} \boxplus_c \sqrt{\mump} \right)^2$.
We have
\begin{equation} \label{subordination_rectangular}
\frac{G_{\nu}(z)}{1-cG_{\nu}(z)} = G_{\mu} \left( z(1-cG_{\nu}(z))^2 - (1-c)(1-cG_{\nu}(z)) \right) \, .
\end{equation}
\end{itemize}
\end{Proposition}

In Theorem \ref{theorem1} below, we are interested in the information-plus-noise model and we control the distance between the spectral measure and the corresponding rectangular free convolution, by bounding the difference between the two terms in (\ref{subordination_rectangular}) evaluated at the average Stieltjes transform.

\subsection{Main results}

Note that in the rest of the paper, we will only consider \underline{real} matrices for ease but our results should generalize to complex matrices adapting the proofs. The only difficulty in the complex case is to adapt the general integration by parts formula (\ref{IPP_general}) which is used several times in this paper, which would lead to heavier computations.\\

Let us define, for $s,t>0$, the distance $d_{s,t}$ on $\PP(\R)$ by
\begin{equation} \label{dst}
d_{s,t}(\mu,\nu) = \sup_{z \in V_{s,t}} \left| G_{\mu}(z) - G_{\nu}(z) \right| \, ,
\end{equation}
where
\begin{equation} \label{Vst}
V_{s,t} = \left\{ z \in \C \ | \ \Imag z > s, \ \left| \frac{\Reel z}{\Imag z} \right| < t \right\} \, .
\end{equation}
As the distance $d$ defined in \cite{BC}, $d_{s,t}$ metrizes weak convergence. Let us mention that for all $\mu,\nu \in \PP(\R)$, we have
\begin{equation} \label{comparison_distances}
d_{s,t}(\mu,\nu) \le \min \left( \dKS(\mu,\nu) , W_1(\mu,\nu) \right) \, ,
\end{equation}
where $\dKS$ and $W_1$ are respectively the Kolmogorov-Smirnov and the $L^1$-Wasserstein distances on $\PP(\R)$. Some key inequalities for the distance between two empirical spectral measures are summarized in Appendix \ref{inequalities_ESM}.

Our first main result is the following.

\begin{Theorem} \label{theorem1}
We assume that $c_n = \frac{n}{p}$ is bounded below and above by two constants in $(0,+\infty)$. Let $c>0$. There exist $s,t>0$ and a constant $c_{s,t}>0$ such that for any random matrix $Y \in \MM_{n,p}(\R)$ with i.i.d. entries satisfying $\Var(Y_{1,1})=1$ and $\E(Y_{1,1}^4) < +\infty$, for any deterministic matrix $M \in \MM_{n,p}(\R)$, and for all $n$ large enough, we have
\begin{multline*}
d_{s,t} \left( \E \mu_{(Y/\sqrt{p}+M)(Y/\sqrt{p}+M)^t} , \left( \sqrt{\mu_{MM^t}} \boxplus_c \sqrt{\mump} \right)^2 \right) \\
\le \ c_{s,t} \left( \E |\mathring{Y_{1,1}}|^3 + \E(\mathring{Y_{1,1}}^4) \right) \left( \frac{1}{\sqrt{n}} + \frac{\Tr(MM^t)^{1/2}}{n} \right) \\
+ \ c_{s,t} \left( |c_n-c| + \frac1n + \frac{\Tr(MM^t)^{1/2}}{n^{5/4}} \right) \, ,
\end{multline*}
where $\mathring{Y}$ is the matrix whose entries are given by $\mathring{Y_{j,k}} = Y_{j,k} - \E(Y_{j,k})$.
\end{Theorem}

This result allows to understand the influence of the deformation in the information-plus-noise model. First, we can observe a decorrelation between the classical term $\frac{1}{\sqrt{n}}$ and the Frobenius norm of the deformation divided by a better power of $n$, namely $\frac{\Tr(MM^t)^{1/2}}{n}$. It is important for us to get this precise estimate since in Section \ref{section_deviations}, we apply Theorem \ref{theorem1} to a matrix $M$ whose Frobenius norm is not bounded but of order $\sqrt{n} \log n$.

Besides, it is interesting to compare Theorem \ref{theorem1} to the Wigner case (cf. \cite[Theorem 2.6]{BC}). Bordenave and Caputo investigated additive deformations and obtained that in this model, the distance between the spectral measure and the corresponding free additive convolution is bounded by $\frac{1}{\sqrt{n}}$. This bound is uniform in the deformation $M$ and it depends on the initial matrix through its moments only. In the case of sample covariance matrices, it would have been surprising if we had obtained a better bound. Table \ref{comparison_bounds} below permits to compare Bordenave and Caputo's results with ours in the Gaussian and the general cases. In addition to this, let us mention that in \cite{CDFF}, the authors were interested in the case of Wigner matrices whose entries have a symmetric distribution satisfying a Poincar\'e inequality, which leads to better bounds than \cite{BC}.

\begin{table}[h]
\begin{tabular}{|c|c|c|}
\hline & Gaussian & Non-Gaussian \\ \hline
\begin{tabular}{c} Wigner \\ matrix \end{tabular}
& \begin{tabular}{c} Deformed GUE matrix \\ $\dfrac1n$ \end{tabular}
& \begin{tabular}{c} Deformed Wigner matrix \\ $\dfrac{1}{\sqrt{n}}$ \end{tabular}
\\ \hline \begin{tabular}{c} Covariance \\ matrix \end{tabular}
& \begin{tabular}{c} Deformed LOE matrix \\ $\dfrac1n + \dfrac{\Tr(MM^t)^{1/2}}{n^{5/4}}$ \end{tabular}
& \begin{tabular}{c} Info-plus-noise matrix \\ $\dfrac{1}{\sqrt{n}} + \dfrac{\Tr(MM^t)^{1/2}}{n}$ \end{tabular} \\ \hline
\end{tabular}
\caption{Bound in the subordination relation (\ref{subordination_square}) or (\ref{subordination_rectangular}) for different matrix models.} \label{comparison_bounds}
\end{table}

Theorem \ref{theorem1} above will be used in the proof of our second main result.

\begin{Theorem} \label{theorem2}
Let $X \in \MM_{n,p}(\R)$ be a random matrix such that $c_n = \frac{n}{p} \to c \in (0,+\infty)$. We assume that $\Var(X_{1,1})=1$ and that there exist $\alpha \in (0,2)$ and $a \in (0,+\infty]$ such that $X_{1,1} \in \Sz_{\alpha}(a)$.\\
Then, the empirical spectral measure $\mu_{XX^t/p}$ satisfies the LDP with speed $n^{1+\alpha/2}$ in $\PP(\R_+)$, governed by the good rate function $J'$ defined by
\begin{displaymath}
J'(\mu) = \left\{ \begin{array}{ll}
\frac{a}{c^{\alpha/2}} m_{\alpha/2}(\nu) & \textrm{ if there exists } \nu \in \PP(\R_+) \textrm{ s.t. } \mu = \left( \sqrt{\nu} \boxplus_c \sqrt{\mump} \right)^2 \\
& \quad \textrm{ and } \nu(\{0\}) \ge \max \left( 0,1-\frac1c \right) \\
+\infty & \textrm{ otherwise}
\end{array} \right.
\end{displaymath}
where $m_p(\mu) = \int_{\R} |x|^p \, d\mu(x)$ denotes the $p$-th moment of a distribution $\mu$.
\end{Theorem}

It is very similar to Bordenave and Caputo's result (see Theorem \ref{LDP_BC}), the main difference being the explicit expression of the rate function in all cases. This is due to the fact that here, we can achieve large deviation explicitly without using a LDP on graphs.\\

The rest of the paper is organized as follows. In Section \ref{section_freeness}, we prove the bound for rectangular free convolution stated in Theorem \ref{theorem1}. In Section \ref{section_deviations}, we prove the large deviation principle in Theorem \ref{theorem2}. In Appendix \ref{appendix_concentration}, we state and prove concentration results used in Sections \ref{section_freeness} and \ref{section_deviations}. Finally, in Appendix \ref{appendix_tools}, we summarize miscellaneous inequalities and identities used throughout the paper.

\section{Asymptotic freeness} \label{section_freeness}

This section is devoted to the proof of Theorem \ref{theorem1}. This theorem is in fact a consequence of the following, as we will see in Section \ref{section_1bis->1}.

\begin{Theorem}[Bound in subordination formula (\ref{subordination_rectangular})] \label{theorem1bis}
We assume that $c_n = \frac{n}{p}$ is bounded below and above by two constants in $(0,+\infty)$. Let $c>0$. There exist $s,t>0$ and a function $f$, bounded on the domain $V_{s,t}$ defined by (\ref{Vst}), such that for any random matrix $Y \in \MM_{n,p}(\R)$ with i.i.d. entries satisfying $\Var(Y_{1,1})=1$ and $\E(Y_{1,1}^4) < +\infty$, for any deterministic matrix $M \in \MM_{n,p}(\R)$, for all $n$ large enough, and for all $z \in V_{s,t}$, we have
\begin{multline*}
\left| \gbar(z) - (1-c\gbar(z)) G_{\mu_{MM^t}}(z(1-c\gbar(z))^2-(1-c)(1-c\gbar(z))) \right| \\
\le \ f(z) \left( \E |\mathring{Y_{1,1}}|^3 + \E(\mathring{Y_{1,1}}^4) \right) \left( \frac{1}{\sqrt{n}} + \frac{\Tr(MM^t)^{1/2}}{n} \right) \\
+ \ f(z) \left( |c_n-c| + \frac1n + \frac{\Tr(MM^t)^{1/2}}{n^{5/4}} \right) \, ,
\end{multline*}
where $g(z) = G_{\mu_{(Y/\sqrt{p}+M)(Y/\sqrt{p}+M)^t}}(z)$ and $\gbar(z) = \E(g(z))$.
\end{Theorem}

The proof of Theorem \ref{theorem1bis} follows the same lines as Bordenave and Caputo's one for the bound in subordination formula (\ref{subordination_square}) for free additive convolution (see \cite[Theorem A.1]{BC}). It consists in two main steps: the Gaussian case and the general case, which we deduce from the Gaussian case by interpolation. However, in the case of sample covariance matrices, the computations are heavier and some majorizations must be finer.

Let us mention that in the Gaussian case, the bound consists only in the last terms (see Proposition \ref{freeness_Gaussian}).

In the proof, we define
\begin{displaymath}
X = \frac{Y}{\sqrt{p}} + M
\end{displaymath}
and we denote by
\begin{displaymath}
S=(zI_n-XX^t)^{-1}
\end{displaymath}
the resolvent of $XX^t$. We consider $s>2$, $t>0$, and along the proof, $s$ can increase and $t$ can decrease. Moreover, $f$ will denote a bounded function on $V_{s,t}$, which can also change from one line to another. In particular, for all $z \in V_{s,t}$ and $x<y$, since we have
\begin{displaymath}
\frac{|z|^x}{|\Imag z|^y} = \frac{1}{|\Imag z|^{y-x}} \left( \left( \frac{\Reel z}{\Imag z} \right)^2 + 1 \right)^{x/2} < \frac{1}{s^{y-x}} (t^2+1)^{x/2} \, ,
\end{displaymath}
we will write
\begin{displaymath}
\frac{|z|^x}{|\Imag z|^y} \le f(z)
\end{displaymath}
as soon as $x<y$.

Before starting the proofs, let us state a lemma we will use in the different steps. $B_{\C}(z,\delta)$ denotes here the ball with centre $z \in \C$ and radius $\delta>0$ for the usual distance in $\C$.

\begin{Lemma} \label{freeness_lemma}
For $\mu \in \PP(\R)$ and $z \in \C$, we define
\begin{displaymath}
\phi_{z,\mu} : (h,\gamma) \mapsto (1-\gamma h) G_{\mu}(z(1-\gamma h)^2-(1-\gamma)(1-\gamma h)) \, .
\end{displaymath}
There exist $s,t>0$ and $l_{s,t}, l'_{s,t} \in (0,1)$ such that
\begin{itemize}
\item for all $\mu \in \PP(\R)$, $z \in V_{s,t}$, and $\gamma>0$, $\phi_{z,\mu}(.,\gamma)$ is Lipschitz on $B_{\C} \left( 0,\frac1s \right)$ with constant $l_{s,t}$,
\item for all $\mu \in \PP(\R)$, $z \in V_{s,t}$, and $h \in B_{\C} \left( 0,\frac1s \right)$, $\phi_{z,\mu}(h,.)$ is Lipschitz on $(0,+\infty)$ with constant $l'_{s,t}$.
\end{itemize}
\end{Lemma}

The proof of this lemma consists in simple computations and is left to the reader. Let us mention however that it relies on the inequality
\begin{equation} \label{ImEta}
|\Imag \eta| > \Imag z \left( \frac{(\sigma-1)(\sigma^2-2)}{\sigma^2(\sigma+1)} - \frac{2t(\sigma+1)}{\sigma^2} \right) - \frac{|1-\gamma|}{\sigma}
\end{equation}
where $\eta = z(1-\gamma h)^2-(1-\gamma)(1-\gamma h)$ and $\sigma = \frac{s}{\gamma}$. We will use it again later.

Furthermore, note that choosing a larger $s$ and a smaller $t$, $l_{s,t}$ and $l'_{s,t}$ can be as close to 0 as wanted.

\subsection{Proof of Theorem \ref{theorem1}} \label{section_1bis->1}

First, let us deduce Theorem \ref{theorem1} from Theorem \ref{theorem1bis}.

\begin{Proof}
We define $\nu = \left( \sqrt{\mu_{MM^t}} \boxplus_c \sqrt{\mump} \right)^2$ and we consider the function $\phi_{z,\mu_{MM^t}}$ defined in Lemma \ref{freeness_lemma}. Subordination formula (\ref{subordination_rectangular}) can be rewritten $\phi_{z,\mu_{MM^t}}(G_{\nu}(z),c) = G_{\nu}(z)$ for all $z \in \C \setminus \R$. Consequently, using Lemma \ref{freeness_lemma}, there exist $s,t>0$ and $l_{s,t} \in (0,1)$ such that for all $z \in V_{s,t}$,
\begin{eqnarray*}
\left| \gbar(z) - G_{\nu}(z) \right|
& \le & \left| \gbar(z) - \phi_{z,\mu_{MM^t}}(\gbar(z),c) \right| \\
& & \quad + \ \left| \phi_{z,\mu_{MM^t}}(\gbar(z),c) - \phi_{z,\mu_{MM^t}}(G_{\nu}(z),c) \right| \\
& \le & \left| \gbar(z) - \phi_{z,\mu_{MM^t}}(\gbar(z),c) \right| + l_{s,t} \left| \gbar(z) - G_{\nu}(z) \right|
\end{eqnarray*}
thus
\begin{displaymath}
\left| \gbar(z) - G_{\nu}(z) \right| \le \frac{1}{1-l_{s,t}} \left| \gbar(z) - \phi_{z,\mu_{MM^t}}(\gbar(z),c) \right| \, .
\end{displaymath}
From Theorem \ref{theorem1bis} in which we majorize $f$ by a constant depending on $s,t$ and from the definition (\ref{dst}) of $d_{s,t}$, we finally get
\begin{multline*}
d_{s,t} \left( \E \mu_{XX^t} , \nu \right) \le c_{s,t} \left( \E |\mathring{Y_{1,1}}|^3 + \E(\mathring{Y_{1,1}}^4) \right) \left( \frac{1}{\sqrt{n}} + \frac{\Tr(MM^t)^{1/2}}{n} \right) \\
+ \ c_{s,t} \left( |c_n-c| + \frac1n + \frac{\Tr(MM^t)^{1/2}}{n^{5/4}} \right) \, .
\end{multline*} \qed
\end{Proof}

\subsection{The Gaussian case}

In this subsection, we assume that $Y_{1,1}$ is a standard Gaussian. Moreover, we will simply denote $g(z)$ and $\gbar(z)$ by $g$ and $\gbar$ respectively (see Theorem \ref{theorem1bis} for their definitions). We will prove the following bound.

\begin{Proposition} \label{freeness_Gaussian}
There exist $s,t>0$ and a function $f$, bounded on $V_{s,t}$, such that for any random matrix $Y \in \MM_{n,p}(\R)$ with i.i.d. standard Gaussian entries, for any deterministic matrix $M \in \MM_{n,p}(\R)$, for all $n$ large enough, and for all $z \in V_{s,t}$, we have
\begin{multline*}
\left| \gbar - (1-c\gbar) G_{\mu_{MM^t}} (z(1-c\gbar)^2 - (1-c)(1-c\gbar)) \right| \\
\le f(z) \left( |c_n-c| + \frac1n + \frac{\Tr(MM^t)^{1/2}}{n^{5/4}} \right) \, .
\end{multline*}
\end{Proposition}

To prove Proposition \ref{freeness_Gaussian}, we will follow and improve some computations by Dumont \emph{et al.}, see \cite[Appendix II]{DHLLN}.

\begin{Lemma}[adaptation from {\cite[Formula (122)]{VLM}}] \label{freeness_Gaussian_lemma1}
Let $Y \in \MM_{n,p}(\R)$ be a random matrix with i.i.d. standard Gaussian entries, let $M \in \MM_{n,p}(\R)$ be a deterministic matrix, and let $z \in \C \setminus \R$. For all integer $n$, we have
\begin{multline} \label{freeness_Gaussian_eq1}
\gbar - \frac1n \Tr(R) = \frac1n \Tr(\Delta R) - \frac{c_n}{n^2} \Tr(\Delta) \Tr(\E(S)R) \\
+ \frac1n \Tr(\Delta' R) - \frac{c_n}{n^2} \Tr(\Delta') \Tr(\E(S)R)
\end{multline}
where
\begin{equation} \label{R}
R = \left( (z(1-c_n\gbar)-1+c_n) I_n - \frac{1}{1-c_n\gbar} MM^t \right)^{-1} \, ,
\end{equation}
\begin{multline} \label{Delta}
\Delta = \frac{1}{p(1-c_n\gbar)} \E \left( \mathring{\overbrace{\Tr(SXM^t)}} \mathring{S} \right) + \frac{c_n z}{1-c_n\gbar} \E( \mathring{g}S ) \\
+ \frac{c_n}{p(1-c_n\gbar)^2} \E \left( \mathring{g} \mathring{\overbrace{\Tr(SXM^t)}} \right) \E(S) \, ,
\end{multline}
and
\begin{multline} \label{Delta'}
\Delta' = \frac{1}{p(1-c_n\gbar)} \E(SXM^tS) + \frac{1}{p(1-c_n\gbar)} \E(zS^2-S) \\
+ \frac{1}{p^2(1-c_n\gbar)^2} \E(\Tr(S^2XM^t)) \E(S) \, .
\end{multline}
\end{Lemma}

In this lemma, we compare $\gbar$ to $\frac1n \Tr(R)$ because, using the notations in Lemma \ref{freeness_lemma}, we have $\frac1n \Tr(R) = \phi_{z,\mu_{MM^t}}(\gbar,c_n)$, so $\frac1n \Tr(R)$ is close to $\phi_{z,\mu_{MM^t}}(\gbar,c)$ by Lemma \ref{freeness_lemma}. That is interesting if we have in mind our goal, which is Proposition \ref{freeness_Gaussian}.

Note that, as \cite[Formula (122)]{VLM}, the proof of Lemma \ref{freeness_Gaussian_lemma1} mainly relies on the Gaussian integration by parts formula (\ref{IPP_Gaussian}), so we do not give it here.

However, we can observe an important difference between Formula (122) in \cite{VLM} and Lemma \ref{freeness_Gaussian_lemma1}, namely the terms in $\Delta'$. In fact, the background here is not exactly the same as in \cite{VLM}. Indeed, Vallet \emph{et al.} consider complex Gaussian entries with independent real and imaginary parts having the same distribution in the matrix $Y$, whereas we consider real Gaussian entries. Consequently, some simplifications do not occur any longer and a new term appears. Behind this phenomenon is the quantity $\zeta = K_{1,1} + 2iK_{1,2} - K_{2,2}$, where $K$ denotes the covariance matrix of the Gaussian vector $(\Reel Y_{1,1}, \Imag Y_{1,1})$. This quantity is equal to 0 in \cite{VLM} and to 1 here, that is why we have an additional term.\\

In the next lemma, we bound the different terms appearing in (\ref{freeness_Gaussian_eq1}). For this, we will use the concentration bounds (\ref{concentration1}) and (\ref{concentration2}) for the terms in $\Delta$ and standard inequalities on traces and resolvents (see Propositions \ref{prop_trace} and \ref{resolvent_covariance}) for the terms in $\Delta'$. Our computations will partially follow those in \cite{VLM}.

\begin{Lemma} \label{freeness_Gaussian_lemma2}
There exist $s,t>0$ and a function $f$, bounded on $V_{s,t}$, such that for all $Y$, $M$, $n$, and $z$ as in Proposition \ref{freeness_Gaussian}, we have
\begin{displaymath}
\left| \gbar - \frac1n \Tr R \right| \le f(z) \left( \frac1n + \frac{\Tr(MM^t)^{1/2}}{n^{5/4}} \right) \, .
\end{displaymath}
\end{Lemma}

This lemma shows that $\frac1n \Tr(R)$ is a deterministic equivalent to the Stieltjes transform $g(z) = \frac1n \Tr(S)$ as soon as $\frac{\Tr(MM^t)^{1/2}}{n^{5/4}}$ tends to 0 as $n \to +\infty$, i.e. when the perturbation $M$ is not too large.

We can compare this result with the bound obtained in \cite[Proposition 6]{VLM}. Two main differences must be highlighted. First, as we mentioned above, the model is not exactly the same. Indeed, we consider real Gaussian entries and not complex Gaussian entries with independent real and imaginary parts, which produces an additional term in $\Delta'$. However, the terms in $\Delta$ are present in both cases, so we can compare the bounds for these terms. Here is the second difference. In \cite{VLM}, the authors assume that $\|M\|$ is uniformly bounded in $n$ and get the bound $\frac{f(z)}{n^2}$. Here, for the terms in $\Delta$, we will get the bound
\begin{displaymath}
f(z) \left( \frac{1}{n^2} + \frac{\Tr(MM^t)^{1/2}}{n^{17/8}} \right) \, .
\end{displaymath}
Moreover, if we use the bound (\ref{concentration2bis}) instead of (\ref{concentration2}) in the proof, and if we observe that $\Tr(MM^t)^{1/2} \le \sqrt{n} \|M\|$, then we get the bound $\frac{f(z)}{n^2} (1+\|M\|)$, which is the same as in \cite{VLM} when $\|M\|$ is uniformly bounded in $n$. Consequently, our bound has two advantages: it is slightly better than the bound in \cite{VLM} and it applies without any assumption on $M$.

\begin{Proof}
First of all, let us remark that
\begin{displaymath}
\left| \frac{1}{1-c_n\gbar} \right| \le f(z)
\end{displaymath}
since $|\gbar| \le \frac{1}{|\Imag z|}$. Besides, we have
\begin{equation} \label{freeness_Gaussian_eq11}
\|R\| \le f(z)
\end{equation}
because on the one hand, $\frac{R}{1-c_n\gbar}$ is a resolvent evaluated at $\eta = z(1-c_n\gbar)^2 - (1-c_n)(1-c_n\gbar)$ so its operator norm is less than $\frac{1}{|\Imag \eta|}$, and on the other hand, we have the inequalities $|1-c_n\gbar| \le 1 + \frac{c_n}{|\Imag z|}$ and (\ref{ImEta}) (we apply the latter with $\sigma = \frac{s}{c_n}$).\\
By proposition \ref{prop_trace} (ii), it follows that
\begin{displaymath}
\left| \frac1n \Tr(\E(S)R) \right| \le \frac{f(z)}{|\Imag z|}
\end{displaymath}
or just
\begin{equation} \label{freeness_Gaussian_eq12}
\left| \frac1n \Tr(\E(S)R) \right| \le f(z) \, .
\end{equation}
Note that more precise bounds can be obtained, see \cite[Appendix E]{VLM}.\\

Next, let us recall that $\Delta$ is defined by
\begin{multline*}
\Delta = \frac{1}{p(1-c_n\gbar)} \E \left( \mathring{\overbrace{\Tr(SXM^t)}} \mathring{S} \right) + \frac{c_n z}{1-c_n\gbar} \E( \mathring{g}S ) \\
+ \frac{c_n}{p(1-c_n\gbar)^2} \E \left( \mathring{g} \mathring{\overbrace{\Tr(SXM^t)}} \right) \E(S)
\end{multline*}
and observe that $\Tr(SXM^t) = \Tr(X^tSM)$.\\

The first term in (\ref{freeness_Gaussian_eq1}) we bound is $\left| \frac{c_n}{n^2} \Tr(\Delta) \Tr(\E(S)R) \right|$. First, using the concentration bounds (\ref{concentration1}), (\ref{concentration2}), and the Cauchy-Schwarz inequality, we get
\begin{eqnarray} \label{freeness_Gaussian_eq13}
\lefteqn{   \left| \frac1n \Tr \left( \frac{1}{p(1-c_n\gbar)} \E \left( \mathring{\overbrace{\Tr(SXM^t)}} \mathring{S} \right) \right) \right|   } \nonumber \\
& = & \left| \frac1n \Tr \left( \frac{1}{p(1-c_n\gbar)} \E \left[ (\Tr(SMX^t)-\E \Tr(SMX^t)) (S-\E S) \right] \right) \right| \nonumber \\
& = & \left| \frac{c_n}{1-c_n\gbar} \E\left[ \frac1n (\Tr(X^tSM)-\E \Tr(X^tSM)) \frac1n (\Tr(S)-\E \Tr(S)) \right] \right| \nonumber \\
& \le & \frac{c_n}{|1-c_n\gbar|} \Var \left( \frac1n \Tr(X^tSM) \right)^{1/2} \Var \left( \frac1n \Tr(S) \right)^{1/2} \nonumber \\
& \le & \frac{c_n}{|1-c_n\gbar|} \left( \frac{9c_n v(z)}{n^{9/4}} \Tr(MM^t) \right)^{1/2} \left( \frac{4c_n u(z)}{n^2} \right)^{1/2} \nonumber \\
& \le & \frac{f(z)}{n^{17/8}} \Tr(MM^t)^{1/2} \, ,
\end{eqnarray}
where $u(z)$ and $v(z)$ are defined in Proposition \ref{concentration1-2}.

Next, using the identity $g = \frac1n\Tr(S)$ and (\ref{concentration1}), we have
\begin{equation} \label{freeness_Gaussian_eq14}
\left| \frac1n \Tr \left( \frac{c_n z}{1-c_n\gbar} \E( \mathring{g}S ) \right) \right| = \left| \frac{c_n z}{1-c_n\gbar} \Var(g) \right| \le \frac{c_n |z|}{|1-c_n\gbar|} \frac{4c_n u(z)}{n^2} \le \frac{f(z)}{n^2}
\end{equation}
where, for the last inequality, we used the definition of $u(z)$ to get
\begin{displaymath}
\frac{|z| u(z)}{|1-c_n\gbar|} \le f(z) \, .
\end{displaymath}

The same arguments also allow to show that
\begin{equation} \label{freeness_Gaussian_eq15}
\left| \frac1n \Tr \left( \frac{c_n}{p(1-c_n\gbar)^2} \E \left( \mathring{g} \mathring{\overbrace{\Tr(SXM^t)}} \right) \E(S) \right) \right| \le \frac{f(z)}{n^{17/8}} \Tr(MM^t)^{1/2} \, .
\end{equation}

Combining inequalities from (\ref{freeness_Gaussian_eq12}) to (\ref{freeness_Gaussian_eq15}) gives
\begin{equation} \label{freeness_Gaussian_eq16}
\left| \frac{c_n}{n^2} \Tr(\Delta) \Tr(\E(S)R) \right| \le f(z) \left( \frac{1}{n^2} + \frac{\Tr(MM^t)^{1/2}}{n^{17/8}} \right) \, .
\end{equation}

Computations are similar for the term $\frac1n \Tr(\Delta R)$, using the additional inequalities (\ref{freeness_Gaussian_eq11}) and $\Tr(RR^*)^{1/2} \le \sqrt{n} \|R\|$ (see Proposition \ref{prop_trace} (iv)). For instance, we have
\begin{eqnarray*}
\lefteqn{   \left| \frac1n \Tr \E \left( \mathring{\overbrace{\Tr(SXM^t)}} \mathring{S}R \right) \right|   } \\
& = & \left| \frac1n \Tr \left( \frac{1}{p(1-c_n\gbar)} \E \left[ (\Tr(SXM^t)-\E \Tr(SXM^t)) (S-\E S)R \right] \right) \right| \\
& \le & \frac{c_n}{|1-c_n\gbar|} \Var \left( \frac1n \Tr(X^tSM) \right)^{1/2} \Var \left( \frac1n \Tr(SR) \right)^{1/2} \\
& \le & \frac{c_n}{|1-c_n\gbar|} \left( \frac{9c_n v(z)}{n^{9/4}} \Tr(MM^t) \right)^{1/2} \left( \frac{4c_n u(z)}{n^{5/2}} \|R\| \Tr(RR^*)^{1/2} \right)^{1/2} \\
& \le & \frac{c_n}{|1-c_n\gbar|} \left( \frac{9c_n v(z)}{n^{9/4}} \Tr(MM^t) \right)^{1/2} \left( \frac{4c_n u(z)}{n^2} f(z)^2 \right)^{1/2} \\
& \le & \frac{f(z)}{n^{17/8}} \Tr(MM^t)^{1/2} \, .
\end{eqnarray*}
Combining with
\begin{displaymath}
\left| \frac1n \Tr \left( \frac{c_n z}{1-c_n\gbar} \E( \mathring{g}SR ) \right) \right| \le \frac{f(z)}{n^2}
\end{displaymath}
and
\begin{displaymath}
\left| \frac1n \Tr \left( c_n \E \left( \mathring{g} \mathring{\overbrace{\Tr(SXM^t)}} \right) \E(S)R \right) \right| \le \frac{f(z)}{n^{17/8}} \Tr(MM^t)^{1/2} \, ,
\end{displaymath}
which have a similar proof, we thus have
\begin{equation} \label{freeness_Gaussian_eq17}
\left| \frac1n \Tr(\Delta R) \right| \le f(z) \left( \frac{1}{n^2} + \frac{\Tr(MM^t)^{1/2}}{n^{17/8}} \right) \, .
\end{equation}

We have bounded the terms in Lemma \ref{freeness_Gaussian_lemma1} in which $\Delta$ appears thanks to the concentration bounds proved in Appendix \ref{appendix_concentration}. We will now consider the terms in which $\Delta'$ appears, in other words the terms not present in \cite{VLM}. To this, we will only use inequalities on traces and resolvents (see Propositions \ref{prop_trace} and \ref{resolvent_covariance}). Let us recall that $\Delta'$ is defined by
\begin{multline*}
\Delta' = \frac{1}{p(1-c_n\gbar)} \E(SXM^tS) + \frac{1}{p(1-c_n\gbar)} \E(zS^2-S) \\
+ \frac{1}{p^2(1-c_n\gbar)^2} \E(\Tr(S^2XM^t)) \E(S) \, .
\end{multline*}

Using inequalities (i)-(iv) in Proposition \ref{prop_trace} and the resolvent identity $SXX^t = zS - I_n$, we get
\begin{eqnarray} \label{freeness_Gaussian_eq18}
\lefteqn{   \left| \Tr(SXM^tS) \right|   } \nonumber \\
& \le & \Tr(SXX^tS^*)^{1/2} \Tr(M^tSS^*M)^{1/2} \nonumber \\
& \le & \left( \frac{n}{|\Imag \overline{z}|} \left( \frac{|z|}{|\Imag z|} + 1 \right) \right)^{1/2} \left( \frac{\sqrt{p}}{|\Imag z|.|\Imag \overline{z}|} \|M\| \Tr(MM^t)^{1/2} \right)^{1/2} \nonumber \\
& \le & n^{3/4} f(z) \Tr(MM^t)^{1/2}
\end{eqnarray}
so
\begin{equation} \label{freeness_Gaussian_eq19}
\left| \frac1n \Tr \left( \frac{1}{p(1-c_n\gbar)} \E(SXM^tS) \right) \right| \le \frac{f(z)}{n^{5/4}} \Tr(MM^t)^{1/2} \, .
\end{equation}

In addition,
\begin{eqnarray} \label{freeness_Gaussian_eq20}
\left| \frac1n \Tr \left( \frac{1}{p(1-c_n\gbar)} \E(zS^2-S) \right) \right| & \le & \frac{1}{np|1-c_n\gbar|} n \left( \frac{|z|}{|\Imag z|^2} + \frac{1}{|\Imag z|} \right) \nonumber \\
& \le & \frac{f(z)}{n}
\end{eqnarray}
and using (\ref{freeness_Gaussian_eq18}) again,
\begin{equation} \label{freeness_Gaussian_eq21}
\left| \frac1n \Tr \left( \frac{1}{p^2(1-c_n\gbar)^2} \E(\Tr(S^2XM^t)) \E(S) \right) \right| \le \frac{f(z)}{n^{5/4}} \Tr(MM^t)^{1/2} \, .
\end{equation}

Consequently, the combination of (\ref{freeness_Gaussian_eq12}), (\ref{freeness_Gaussian_eq19}), (\ref{freeness_Gaussian_eq20}), and (\ref{freeness_Gaussian_eq21}) gives
\begin{equation} \label{freeness_Gaussian_eq22}
\left| \frac{c_n}{n^2} \Tr(\Delta') \Tr(\E(S)R) \right| \le f(z) \left( \frac1n + \frac{\Tr(MM^t)^{1/2}}{n^{5/4}} \right) \, .
\end{equation}

By very similar calculations, we get
\begin{equation} \label{freeness_Gaussian_eq23}
\left| \frac1n \Tr(\Delta'R) \right| \le f(z) \left( \frac1n + \frac{\Tr(MM^t)^{1/2}}{n^{5/4}} \right) \, .
\end{equation}

Finally, combining relation (\ref{freeness_Gaussian_eq1}) with inequalities (\ref{freeness_Gaussian_eq16}), (\ref{freeness_Gaussian_eq17}), (\ref{freeness_Gaussian_eq22}), and (\ref{freeness_Gaussian_eq23}), we get
\begin{displaymath}
\left| \gbar - \frac1n \Tr(R) \right| \le f(z) \left( \frac1n + \frac{\Tr(MM^t)^{1/2}}{n^{5/4}} \right) \ .
\end{displaymath} \qed
\end{Proof}

Finally, the Gaussian case (Proposition \ref{freeness_Gaussian}) follows from Lemma \ref{freeness_Gaussian_lemma2} and the second part of Lemma \ref{freeness_lemma}, since we have
\begin{displaymath}
\frac1n \Tr(R) = (1-c_n\gbar) G_{\mu_{MM^t}}(z(1-c_n\gbar)^2 - (1-c_n)(1-c_n\gbar)) \, .
\end{displaymath}

\subsection{The general case}

We now only assume that $\Var(Y_{1,1})=1$ and that $\E(Y_{1,1}^4) < +\infty$. Let $\Yhat \in \MM_{n,p}(\R)$ be an independent random matrix such that the $\Yhat_{j,k}$'s are i.i.d. standard Gaussians, we define $\Xhat = \frac{\Yhat}{\sqrt{p}} + M$ and for all $u \in [0,1]$, we define $Y(u) = \sqrt{u}Y + \sqrt{1-u}\Yhat$, $X(u) = \frac{Y(u)}{\sqrt{p}} + M$, and $S(u) = (zI_n - X(u)X(u)^t)^{-1}$. We have the following, which will allow us to bring back the general case to the Gaussian case.

\begin{Proposition} \label{freeness_general_step1}
There exist $s,t>0$ and a function $f$, bounded on $V_{s,t}$, such that for any random matrix $Y \in \MM_{n,p}(\R)$ with i.i.d. entries satisfying $\Var(Y_{1,1})=1$, $\E(Y_{1,1}^4) < +\infty$, \underline{and $\E(Y_{1,1})=0$}, for any deterministic matrix $M \in \MM_{n,p}(\R)$, for all $n$ large enough, and for all $z \in V_{s,t}$, we have
\begin{multline*}
\left| \E G_{\mu_{XX^t}}(z) - \E G_{\mu_{\Xhat\Xhat^t}}(z) \right| \\
\le f(z) \left( \E |Y_{1,1}|^3 + \E(Y_{1,1}^4) \right) \left( \frac{1}{\sqrt{n}} + \frac{\Tr(MM^t)^{1/2}}{n} \right) \, .
\end{multline*}
\end{Proposition}

\begin{Proof}
The proof consists in four main steps. After developing $\E G_{\mu_{XX^t}}(z) - \E G_{\mu_{\Xhat\Xhat^t}}(z)$, we use integration by parts formulas (see Lemma \ref{IPP_lemma}). Then, we respectively focus on bounds for the main terms and the rests in these integrations by parts.\\

\noindent \boxed{\textbf{First step: Development of $\E G_{\mu_{XX^t}}(z) - \E G_{\mu_{\Xhat\Xhat^t}}(z)$.}}\\
Let $u \in [0,1]$ and $h \in [-u,1-u]$. Proposition \ref{resolvent_covariance} (ii), applied to $A=X(u)$ and $B = \frac{1}{\sqrt{p}} (Y(u+h)-Y(u))$ gives
\begin{eqnarray*}
S(u+h)-S(u) & = & S(u+h) \left( X(u) \left( \frac{Y(u+h)-Y(u)}{\sqrt{p}} \right)^t \right. \\
& & \quad + \frac{Y(u+h)-Y(u)}{\sqrt{p}} X(u)^t \\
& & \quad \left. + \frac{Y(u+h)-Y(u)}{\sqrt{p}} \left( \frac{Y(u+h)-Y(u)}{\sqrt{p}} \right)^t \right) S(u) \, .
\end{eqnarray*}
Dividing by $h$ and taking $h \to 0$, we get for all $u \in [0,1]$,
\begin{eqnarray*}
S'(u) & = & S(u) \left( X(u) \frac{Y'(u)^t}{\sqrt{p}} + \frac{Y'(u)}{\sqrt{p}} X(u)^t \right) S(u) \\
& = & \frac{1}{\sqrt{p}} S(u) \left( \left( \frac{\sqrt{u}Y}{\sqrt{p}} + \frac{\sqrt{1-u}\Yhat}{\sqrt{p}} + M \right) \left( \frac{Y^t}{2\sqrt{u}} - \frac{\Yhat^t}{2\sqrt{1-u}} \right) \right. \\
& & \quad \left. + \left( \frac{Y}{2\sqrt{u}} - \frac{\Yhat}{2\sqrt{1-u}} \right) \left( \frac{\sqrt{u}Y^t}{\sqrt{p}} + \frac{\sqrt{1-u}\Yhat^t}{\sqrt{p}} + M^t \right) \right) S(u) \, .
\end{eqnarray*}
Thus we can rewrite
\begin{eqnarray*}
\lefteqn{  G_{\mu_{XX^t}}(z) - G_{\mu_{\Xhat\Xhat^t}}(z)  } \\
& = & \frac1n \Tr S(1) - \frac1n \Tr S(0) \\
& = & \frac1n \int_0^1 \Tr S'(u) \, du \\
& = & \frac{1}{2n\sqrt{p}} \int_0^1 \Tr \left[ S(u)^2 \left[ 2\frac{YY^t}{\sqrt{p}} - 2\frac{\Yhat\Yhat^t}{\sqrt{p}} \right. \right. \\
& & \quad + \left( \sqrt{\frac{1-u}{u}} - \sqrt{\frac{u}{1-u}} \right) \frac{Y\Yhat^t}{\sqrt{p}} + \left( \sqrt{\frac{1-u}{u}} - \sqrt{\frac{u}{1-u}} \right) \frac{\Yhat Y^t}{\sqrt{p}} \\
& & \quad \left. \left. + \frac{MY^t}{\sqrt{u}} - \frac{M\Yhat^t}{\sqrt{1-u}} + \frac{YM^t}{\sqrt{u}} - \frac{\Yhat M^t}{\sqrt{1-u}} \right] \right] \, du \, .
\end{eqnarray*}
Denoting by
\begin{eqnarray*}
(1) & = & \Tr S(u)^2 \left[ \frac{YY^t}{\sqrt{p}} - \sqrt{\frac{u}{1-u}} \frac{Y\Yhat^t}{\sqrt{p}} \right] \\
& = & \sum_{\substack{1 \le j,k \le n \\ 1 \le l \le p}} \frac{1}{\sqrt{p}} \left[ S(u)^2_{j,k} Y_{k,l} Y_{j,l} - \sqrt{\frac{u}{1-u}} S(u)^2_{j,k} Y_{k,l} \Yhat_{j,l} \right] \, ,
\end{eqnarray*}
\begin{displaymath}
(2) = \sum_{\substack{1 \le j,k \le n \\ 1 \le l \le p}} \frac{1}{\sqrt{p}} \left[ S(u)^2_{j,k} Y_{k,l} Y_{j,l} - \sqrt{\frac{u}{1-u}} S(u)^2_{j,k} \Yhat_{k,l} Y_{j,l} \right] \, ,
\end{displaymath}
\begin{displaymath}
(3) = \sum_{\substack{1 \le j,k \le n \\ 1 \le l \le p}} \frac{1}{\sqrt{p}} \left[ \sqrt{\frac{1-u}{u}} S(u)^2_{j,k} Y_{k,l} \Yhat_{j,l} - S(u)^2_{j,k} \Yhat_{k,l} \Yhat_{j,l} \right] \, ,
\end{displaymath}
\begin{displaymath}
(4) = \sum_{\substack{1 \le j,k \le n \\ 1 \le l \le p}} \frac{1}{\sqrt{p}} \left[ \sqrt{\frac{1-u}{u}} S(u)^2_{j,k} \Yhat_{k,l} Y_{j,l} - S(u)^2_{j,k} \Yhat_{k,l} \Yhat_{j,l} \right] \ ,
\end{displaymath}
\begin{displaymath}
(5) = \sum_{\substack{1 \le j,k \le n \\ 1 \le l \le p}} \frac{1}{\sqrt{u}} S(u)^2_{j,k} M_{k,l} Y_{j,l} - \frac{1}{\sqrt{1-u}} S(u)^2_{j,k} M_{k,l} \Yhat_{j,l} \, ,
\end{displaymath}
and
\begin{displaymath}
(6) = \sum_{\substack{1 \le j,k \le n \\ 1 \le l \le p}} \frac{1}{\sqrt{u}} S(u)^2_{j,k} Y_{k,l} M_{j,l} - \frac{1}{\sqrt{1-u}} S(u)^2_{j,k} \Yhat_{k,l} M_{j,l} \, ,
\end{displaymath}
where $S(u)^2_{j,k}$ must be read $(S(u)^2)_{j,k}$, we finally rewrite
\begin{equation} \label{freeness_general_eq1}
\E G_{\mu_{XX^t}}(z) - \E G_{\mu_{\Xhat\Xhat^t}}(z) = \frac{1}{2n\sqrt{p}} \int_0^1 \E[(1)+(2)+(3)+(4)+(5)+(6)] \, du \, .
\end{equation}\\

\noindent \boxed{\textbf{Second step: Integrations by parts.}}\\
Let us recall the formulas we will use below.

\begin{Lemma}[see {\cite[Formulas (2.1.39) and (18.1.19)]{PS}}] \label{IPP_lemma}
\begin{itemize}
\item[(i)] Let a function $F \in \mathcal{C}^1(\R,\R)$ and $\xi$ a random variable with distribution $\NN(0,\sigma^2)$. If $\E |F'(\xi)| < +\infty$, then
\begin{equation} \label{IPP_Gaussian}
\E \left( F(\xi)\xi \right) = \sigma^2 \E \left( F'(\xi) \right) \, .
\end{equation}
\item[(ii)] More generally, let $p$ be an integer, a function $F \in \mathcal{C}^{p+1}(\R,\R)$, and a real random variable $\xi$. If $\E |\xi|^{p+2} < +\infty$ and the derivatives $F',\ldots,F^{(p+1)}$ are bounded on $\R$, then
\begin{equation} \label{IPP_general}
\E \left( \xi F(\xi) \right) = \sum_{j=0}^p \frac{\kappa_{j+1}}{j!} \E(F^{(j)}(\xi)) + \varepsilon_p
\end{equation}
where the $\kappa_{j+1}$'s are the cumulants of the distribution of $\xi$ and
\begin{displaymath}
|\varepsilon_p| \le C_p \E |\xi|^{p+2} . \|F^{(p+1)}\|_{\infty} \ , \qquad C_p \le \frac{1+(3+2p)^{p+2}}{(p+1)!} \, .
\end{displaymath}
\end{itemize}
\end{Lemma}

We will apply the Gaussian (\ref{IPP_Gaussian}) or the general (\ref{IPP_general}) integration by parts formula for all $j,k,l$ in order to decompose $\E[(1)+(2)+(3)+(4)+(5)+(6)]$ as a sum of terms that we can bound.

Note a first crucial point here. As we want to apply Theorem \ref{theorem1} to the matrices $Y$ and $C$ in Section \ref{section_deviations} in order to obtain (\ref{exp_equiv_eq7}), it will not be sufficient to use the integration by parts formula up to order 2, that is why we will be interested in terms of order 3 in this formula.

From now, $D_{a,b}$ denotes the derivation with respect to $Y_{a,b}$.\\

Let $u \in [0,1]$, $j,k \in \llbracket 1,n \rrbracket$, and $l \in \llbracket 1,p \rrbracket$. We denote by $F_1$ and $G_1$ the functions defined by $F_1(Y_{j,l}) = Y_{k,l}S(u)^2_{j,k}$ and $G_1(\Yhat_{j,l}) = Y_{k,l}S(u)^2_{j,k}$. We have
\begin{displaymath}
F_1'(Y_{j,l}) = \frac{2\sqrt{u}}{\sqrt{p}} Y_{k,l}S(u)_{j,k}.D_{j,l}S(u)_{j,k} + \delta_{j,k} S(u)^2_{k,k} \, ,
\end{displaymath}
\begin{multline*}
F_1''(Y_{j,l}) = \frac{2u}{p} Y_{k,l} \left( (D_{j,l}S(u)_{j,k})^2 + S(u)_{j,k}.D^2_{j,l}S(u)_{j,k} \right) \\
+ \frac{4\sqrt{u}}{\sqrt{p}} \delta_{j,k} S(u)_{k,k}.D_{k,l}S(u)_{k,k} \, ,
\end{multline*}
\begin{multline*}
F_1^{(3)}(Y_{j,l}) = \frac{2u^{3/2}}{p^{3/2}} Y_{k,l} \left( 3D_{j,l}S(u)_{j,k}.D^2_{j,l}S(u)_{j,k} + S(u)_{j,k}.D^3_{j,l}S(u)_{j,k} \right) \\
+ \, \frac{6u}{p} \delta_{j,k} \left( (D_{k,l}S(u)_{k,k})^2 + S(u)_{k,k}.D^2_{k,l}S(u)_{k,k} \right) \, ,
\end{multline*}
and
\begin{displaymath}
G_1'(\Yhat_{j,l}) = \frac{2\sqrt{1-u}}{\sqrt{p}} Y_{k,l}S(u)_{j,k}.D_{j,l}S(u)_{j,k} \, .
\end{displaymath}
Applied conditionally to the variables $\{ Y_{a,b}, \ 1 \le a \le n, \ 1 \le b \le p \} \cup \{ \Yhat_{a,b}, \ (a,b) \neq (j,l) \}$, (\ref{IPP_Gaussian}) gives
\begin{displaymath}
\widehat{\E}_{j,l} (S(u)^2_{j,k} Y_{k,l} \Yhat_{j,l}) = \Var(\Yhat_{j,l}) \widehat{\E}_{j,l}(G_1'(\Yhat_{j,l})) \, ,
\end{displaymath}
where $\widehat{\E}_{j,l}$ denotes the associated conditional expectation.
Similarly, from (\ref{IPP_general}), we have
\begin{displaymath}
\E_{j,l} (S(u)^2_{j,k} Y_{k,l} Y_{j,l}) = \Var(Y_{j,l}) \E_{j,l}(F_1'(Y_{j,l})) + \frac{\kappa_3(Y_{j,l})}{2} \E_{j,l}(F_1''(Y_{j,l})) + \varepsilon_{1,j,k,l} \, ,
\end{displaymath}
where $\E_{j,l}$ denotes the expectation conditionally to the variables $\{ Y_{a,b}, \ (a,b) \neq (j,l) \} \cup \{ \Yhat_{a,b}, \ 1 \le a \le n, \ 1 \le b \le p \}$.\\
Taking the expectation, we thus have
\begin{eqnarray*}
\lefteqn{   \E \left[ S(u)^2_{j,k} Y_{k,l} Y_{j,l} - \sqrt{\frac{u}{1-u}} S(u)^2_{j,k} Y_{k,l} \Yhat_{j,l} \right]   } \\
& = & \Var(Y_{j,l}) \E(F_1'(Y_{j,l})) + \frac{\kappa_3(Y_{j,l})}{2} \E(F_1''(Y_{j,l})) + \E(\varepsilon_{1,j,k,l}) \\
& & \quad - \sqrt{\frac{u}{1-u}} \Var(\Yhat_{j,l}) \E(G_1'(\Yhat_{j,l})) \\
& = & \delta_{j,k} \E(S(u)^2_{k,k}) + \frac{\kappa_3(Y_{j,l})}{2} \E(F_1''(Y_{j,l})) + \E(\varepsilon_{1,j,k,l})
\end{eqnarray*}
with
\begin{displaymath}
|\varepsilon_{1,j,k,l}| \le \frac{1+7^4}{6} \E(Y_{1,1}^4) . \|F_1^{(3)}\|_{\infty} \, .
\end{displaymath}
Dividing by $\sqrt{p}$ and summing on $j,k,l$, we thus have
\begin{eqnarray} \label{freeness_general_eq2}
\E(1) & = & \underbrace{  \sqrt{p} \E (\Tr S(u)^2)  }_{(1.1)}  +  \underbrace{  \frac{\kappa_3(Y_{1,1})u}{p^{3/2}} \sum_{j,k,l} \E (Y_{k,l}(D_{j,l}S(u)_{j,k})^2)  }_{(1.2)} \nonumber \\
& & \quad +  \underbrace{  \frac{\kappa_3(Y_{1,1})u}{p^{3/2}} \sum_{j,k,l} \E (Y_{k,l}S(u)_{j,k}.D^2_{j,l}S(u)_{j,k})  }_{(1.3)} \nonumber \\
& & \quad +  \underbrace{  \frac{2\kappa_3(Y_{1,1})\sqrt{u}}{p} \sum_{k,l} \E (S(u)_{k,k}.D_{k,l}S(u)_{k,k})  }_{(1.4)} \nonumber \\
& & \quad +  \frac{1}{\sqrt{p}} \sum_{j,k,l} \E(\varepsilon_{1,j,k,l}) \, .
\end{eqnarray}

Since $S(u)^t = S(u)$, we also have
\begin{equation} \label{freeness_general_eq3}
\E(2) = (1.1) + (1.2) + (1.3) + (1.4) + \frac{1}{\sqrt{p}} \sum_{j,k,l} \E(\varepsilon_{2,j,k,l})
\end{equation}
with
\begin{displaymath}
|\varepsilon_{2,j,k,l}| \le \frac{1+7^4}{6} \E(Y_{1,1}^4) . \|F_1^{(3)}\|_{\infty} \, .
\end{displaymath}

Similarly, considering $F_3(Y_{k,l}) = \Yhat_{j,l} S(u)^2_{j,k}$, we get
\begin{eqnarray} \label{freeness_general_eq4}
\E(3) & = & \underbrace{  -\sqrt{p} \E (\Tr S(u)^2)  }_{(3.1)}  +  \underbrace{  \frac{\kappa_3(Y_{1,1})\sqrt{u(1-u)}}{p^{3/2}} \sum_{j,k,l} \E (\Yhat_{j,l}(D_{k,l}S(u)_{j,k})^2)  }_{(3.2)} \nonumber \\
& & \quad + \underbrace{  \frac{\kappa_3(Y_{1,1})\sqrt{u(1-u)}}{p^{3/2}} \sum_{j,k,l} \E (\Yhat_{j,l}S(u)_{j,k}.D^2_{k,l}S(u)_{j,k})  }_{(3.3)} \nonumber \\
& & \quad +  \frac{1}{\sqrt{p}} \sum_{j,k,l} \E(\varepsilon_{3,j,k,l}) \, ,
\end{eqnarray}
\begin{equation} \label{freeness_general_eq5}
\E(4) = (3.1) + (3.2) + (3.3) + \frac{1}{\sqrt{p}} \sum_{j,k,l} \E(\varepsilon_{4,j,k,l})
\end{equation}
with
\begin{displaymath}
|\varepsilon_{3,j,k,l}| \le \frac{1+7^4}{6} \E(Y_{1,1}^4) . \|F_3^{(3)}\|_{\infty} \quad \textrm{and} \quad |\varepsilon_{4,j,k,l}| \le \frac{1+7^4}{6} \E(Y_{1,1}^4) . \|F_3^{(3)}\|_{\infty}
\end{displaymath}
for all $j,k,l$, and considering $F_5(Y_{j,l}) = M_{k,l} S(u)^2_{j,k}$,
\begin{eqnarray} \label{freeness_general_eq6}
\E(5) & = & \underbrace{  \frac{\kappa_3(Y_{1,1})\sqrt{u}}{p} \sum_{j,k,l} \E (M_{k,l}(D_{j,l}S(u)_{j,k})^2)  }_{(5.1)} \nonumber \\
& & \quad  +  \underbrace{  \frac{\kappa_3(Y_{1,1})\sqrt{u}}{p} \sum_{j,k,l} \E (M_{k,l}S(u)_{j,k}.D^2_{j,l}S(u)_{j,k})  }_{(5.2)} \nonumber \\
& & \quad +  \sum_{j,k,l} \E(\varepsilon_{5,j,k,l}) \, ,
\end{eqnarray}
\begin{equation} \label{freeness_general_eq7}
\E(6) = (5.1) + (5.2) + \sum_{j,k,l} \E(\varepsilon_{6,j,k,l})
\end{equation}
with
\begin{displaymath}
|\varepsilon_{5,j,k,l}| \le \frac{1+7^4}{6} \E(Y_{1,1}^4) . \|F_5^{(3)}\|_{\infty} \quad \textrm{and} \quad |\varepsilon_{6,j,k,l}| \le \frac{1+7^4}{6} \E(Y_{1,1}^4) . \|F_5^{(3)}\|_{\infty}
\end{displaymath}
for all $j,k,l$.

We have thus rewritten
\begin{multline*}
\E G_{\mu_{XX^t}}(z) - \E G_{\mu_{\Xhat\Xhat^t}}(z) \\
= \frac{1}{n\sqrt{p}} \int_0^1 \E[(1.2)+(1.3)+(1.4)+(3.2)+(3.3)+(5.1)+(5.2)] \, du \\
+ \frac{1}{2n\sqrt{p}} \int_0^1 \E \left[ \sum_{j,k,l} \left( \frac{1}{\sqrt{p}} \sum_{i=1}^4 \varepsilon_{i,j,k,l} + \sum_{i=5}^6 \varepsilon_{i,j,k,l} \right) \right] \, du \, .
\end{multline*}\\

\noindent \boxed{\textbf{Third step: Bounds for the main terms.}}\\
We will develop the different terms in this expression with the differentiation formulas in Proposition \ref{resolvent_covariance} (vii), and bound them thanks to inequalities on traces and resolvents (see Propositions \ref{prop_trace} and \ref{resolvent_covariance} again).

Because some computations are very similar, we will be interested in the terms (1.2), (1.3), and (1.4) only.

Note that in order to simplify the notations, from now, we will denote $S$ and $X$ for $S(u)$ and $X(u)$.\\

Let us start with the term (1.2). Using (\ref{DabSjk}) and (\ref{D2abSjk}), we have
\begin{eqnarray*}
\lefteqn{   \sum_{j,k,l} \E (Y_{k,l}(D_{j,l}S(u)_{j,k})^2)   } \\
& = & \sum_{j,k,l} \E \left[ Y_{k,l} (SX)_{j,l}^2 S_{j,k}^2 + 2 Y_{k,l} (SX)_{j,l} S_{j,k} S_{j,j} (SX)_{k,l} + Y_{k,l} S_{j,j}^2 (SX)_{k,l}^2 \right] \\
& = & \E \Bigg[ \Tr( Y(X^tS)^{\circ 2} S^{\circ 2} ) + 2 \Tr( (Y \circ SX)X^tS\diag(S)S ) \\
& & \quad + \, \sum_j S_{j,j}^2 . \sum_{k,l} Y_{k,l} (SX)_{k,l}^2 \Bigg]
\end{eqnarray*}
where $\circ$ is the Hadamard product (see Appendix \ref{traces_norms}) and $S^{\circ 2}$ denotes $S \circ S$. Note that it is crucial here to rewrite precisely the terms with the Hadamard product and then to bound the traces rather than bound directly the entries. Indeed, it allows us to get better powers of $n$ in the bound, which is crucial if we have in mind the large deviations in Section \ref{section_deviations}.

Using Propositions \ref{prop_trace}, \ref{resolvent_covariance}, and the Cauchy-Schwarz inequality in $\C^{np}$ and denoting by $y$ a square root of $z$, we have
\begin{eqnarray*}
|\Tr( Y(X^tS)^{\circ 2} S^{\circ 2} )| & \le & \sqrt{p} \|X^tS\|^2.\|S\|^2 \Tr(YY^t)^{1/2} \\
& \le & \frac{\sqrt{p}}{|\Imag y|^2 |\Imag z|^2} \Tr(YY^t)^{1/2} \ ,
\end{eqnarray*}
\begin{eqnarray*}
\lefteqn{   |\Tr( (Y \circ SX)X^tS\diag(S)S )|   } \\
& \le & \sqrt{p} \|X^tS\|.\|\diag(S)\|.\|S\| \Tr((Y \circ SX)(Y \circ SX)^*)^{1/2} \\
& \le & \frac{\sqrt{p}}{|\Imag y|^2 |\Imag z|^2} \Tr(YY^t)^{1/2} \, ,
\end{eqnarray*}
and
\begin{eqnarray*}
\left| \sum_j S_{j,j}^2 . \sum_{k,l} Y_{k,l} (SX)_{k,l}^2 \right|
& \le & \frac{n}{|\Imag z|^2 |\Imag y|} \sum_{k,l} |Y_{k,l} (SX)_{k,l}| \\
& \le & \frac{n}{|\Imag z|^2 |\Imag y|} \left( \sum_{k,l} Y_{k,l}^2 \right)^{1/2} \left( \sum_{k,l} |(SX)_{k,l}|^2 \right)^{1/2} \\
& = & \frac{n \Tr(YY^t)^{1/2}}{|\Imag z|^2 |\Imag y|} \Tr(SXX^tS^*)^{1/2} \\
& \le & \frac{n \Tr(YY^t)^{1/2}}{|\Imag z|^2 |\Imag y|} \left( \frac{n|z|}{|\Imag z|^2} + \frac{n}{|\Imag z|} \right)^{1/2} \, .
\end{eqnarray*}
Using also the bound (\ref{Imy}), there exists a function $f$, bounded on $V_{s,t}$, independent from $Y$, $M$, and $n$, such that for all $z \in V_{s,t}$, we have
\begin{displaymath}
|(1.2)| \le |\kappa_3(Y_{1,1})| f(z) \E (\Tr(YY^t)^{1/2}) \, .
\end{displaymath}
But for a centred random variable, the third cumulant equals the third moment, so this inequality can be rewritten
\begin{displaymath}
|(1.2)| \le \E|Y_{1,1}|^3 f(z) \E (\Tr(YY^t)^{1/2}) \, .
\end{displaymath}

We adopt the same strategy for the term (1.3). We have
\begin{eqnarray*}
\lefteqn{   \sum_{j,k,l} \E (Y_{k,l}S(u)_{j,k}.D^2_{j,l}S(u)_{j,k})   } \\
& = & \sum_{j,k,l} \E[ Y_{k,l} S_{j,k} . 2(S_{j,j} S_{j,k} + (SX)_{j,l}^2 S_{j,k}  \\
& & \quad + \, S_{j,j} (X^tSX)_{l,l} S_{j,k} + 2 S_{j,j} (SX)_{j,l} (SX)_{k,l}) ] \\
& = & 2 \E \left[ \sum_{j,l} S_{j,j} (S^{\circ 2}Y)_{j,l} + \Tr( Y(X^tS)^{\circ 2}S^{\circ 2} ) \right. \\
& & \quad \left. + \sum_{j,l} S_{j,j} (S^{\circ 2}Y)_{j,l} (X^tSX)_{l,l} + 2 \Tr( (Y \circ SX)X^tS\diag(S)S ) \right]
\end{eqnarray*}
so, using the previous bounds, and also
\begin{eqnarray*}
\left| \sum_{j,l} S_{j,j} (S^{\circ 2}Y)_{j,l} \right|
& \le & \frac{1}{|\Imag z|} \sqrt{np} \left( \sum_{j,l} |(S^{\circ 2}Y)_{j,l}|^2 \right)^{1/2} \\
& = & \frac{\sqrt{np}}{|\Imag z|} \Tr( S^{\circ 2}YY^t(S^{\circ 2})^* )^{1/2} \\
& \le & \frac{\sqrt{np}}{|\Imag z|} \left( \sqrt{n} \frac{\|Y\|}{|\Imag z|^4} \Tr(YY^t)^{1/2} \right)^{1/2} \\
& \le & \frac{n^{3/4} p^{1/2}}{|\Imag z|^3} \Tr(YY^t)^{1/2}
\end{eqnarray*}
and
\begin{eqnarray*}
\left| \sum_{j,l} S_{j,j} (S^{\circ 2}Y)_{j,l} (X^tSX)_{l,l} \right|
& \le & \frac{1}{|\Imag z|} \left( 1 + \left| \frac{z}{\Imag z} \right| \right) \sqrt{np} \left( \sum_{j,l} |(S^{\circ 2}Y)_{j,l}|^2 \right)^{1/2} \\
& \le & \frac{n^{3/4} p^{1/2}}{|\Imag z|^3} \left( 1 + \left| \frac{z}{\Imag z} \right| \right) \Tr(YY^t)^{1/2} \, ,
\end{eqnarray*}
the same arguments as above lead to
\begin{displaymath}
|(1.3)| \le \frac{\E|Y_{1,1}|^3}{n^{1/4}} f(z) \E (\Tr(YY^t)^{1/2}) \, .
\end{displaymath}

Besides, we have
\begin{displaymath}
\sum_{k,l} \E (S(u)_{k,k}.D_{k,l}S(u)_{k,k}) = \sum_{k,l} \E[ S_{k,k} . 2 S_{k,k} (SX)_{k,l}]
\end{displaymath}
and
\begin{eqnarray*}
\left| \sum_{k,l} S_{k,k}^2 (SX)_{k,l} \right|
& \le & \frac{1}{|\Imag z|^2} \sqrt{np} \left( \sum_{k,l} |(SX)_{k,l}|^2 \right)^{1/2} \\
& = & \frac{\sqrt{np}}{|\Imag z|^2} \Tr(SXX^tS^*)^{1/2} \\
& \le & \frac{\sqrt{np}}{|\Imag z|^2} \left( \frac{n|z|}{|\Imag z|^2} + \frac{n}{|\Imag z|} \right)^{1/2}
\end{eqnarray*}
thus we get
\begin{displaymath}
|(1.4)| \le \E|Y_{1,1}|^3 f(z) \sqrt{n} \, .
\end{displaymath}

We finally have
\begin{equation}
|(1.2)+(1.3)+(1.4)| \le \E|Y_{1,1}|^3 f(z) \left( \E (\Tr(YY^t)^{1/2}) + \sqrt{n} \right) \, .
\end{equation}

Very similar computations allow to show that
\begin{displaymath}
|(3.2)+(3.3)| \le \E|Y_{1,1}|^3 f(z) \E (\Tr(\Yhat\Yhat^t)^{1/2})
\end{displaymath}
and
\begin{equation}
|(5.1)+(5.2)| \le \E|Y_{1,1}|^3 f(z) \sqrt{n} \Tr(MM^t)^{1/2} \, .
\end{equation}
If we remember that $Y_{1,1}$ and $\Yhat_{1,1}$ have mean zero and variance 1, we have $\E (\Tr(YY^t))^{1/2} \le \sqrt{np}$ and $\E (\Tr(\Yhat\Yhat^t))^{1/2} \le \sqrt{np}$ by Jensen's inequality. Finally, we can write
\begin{multline} \label{freeness_general_eq10}
|(1.2)+(1.3)+(1.4)+(3.2)+(3.3)+(5.1)+(5.2)| \\
\le \E|Y_{1,1}|^3 f(z) \left( n + \sqrt{n} \Tr(MM^t)^{1/2} \right) \, .
\end{multline}\\

\noindent \boxed{\textbf{Fourth step: Bounds for the rests.}}\\
The only thing to be left is to bound the rests appeared in the integration by parts formulas.
We recall that for all $j,k \in \llbracket1,n \rrbracket$, $l \in \llbracket 1,p \rrbracket$, we have
\begin{displaymath}
|\varepsilon_{1,j,k,l}| \le \frac{1+7^4}{6} \E(Y_{1,1}^4) . \|F_1^{(3)}\|_{\infty} \, .
\end{displaymath}
Using the expression of $F_1^{(3)}(Y_{j,l})$, differentiation formulas (\ref{DabSjk}), (\ref{D2abSjk}), (\ref{D3abSjk}), and inequalities (iv)-(vi) in Proposition \ref{resolvent_covariance}, there exists a function $f$, independent from $Y,M,n,j,k,l$, bounded on $V_{s,t}$, such that
\begin{displaymath}
|\varepsilon_{1,j,k,l}| \le f(z) \E(Y_{1,1}^4) \left( \frac{u^{3/2}}{p^{3/2}} |Y_{k,l}| + \frac{u}{p} \delta_{j,k} \right) \, .
\end{displaymath}
So, using the Cauchy-Schwarz inequality in $\R^{np}$, we have
\begin{eqnarray} \label{freeness_general_eq11}
\frac{1}{\sqrt{p}} \left| \sum_{j,k,l} \E(\varepsilon_{1,j,k,l}) \right|
& \le & f(z) \E(Y_{1,1}^4) \left( \frac1p \sum_{k,l} \E |Y_{k,l}| + \sqrt{n} \right) \nonumber \\
& \le & f(z) \E(Y_{1,1}^4) \left( \E (\Tr(YY^t)^{1/2}) + \sqrt{n} \right) \nonumber \\
& \le & f(z) \E(Y_{1,1}^4) n \, .
\end{eqnarray}
The same bound holds for $\frac{1}{\sqrt{p}} \left| \sum_{j,k,l} \E(\varepsilon_{2,j,k,l}) \right|$. Similarly, we get
\begin{equation} \label{freeness_general_eq12}
\frac{1}{\sqrt{p}} \left| \sum_{j,k,l} \E(\varepsilon_{3,j,k,l}) + \E(\varepsilon_{4,j,k,l}) \right| \le f(z) \E(Y_{1,1}^4) n
\end{equation}
and
\begin{equation} \label{freeness_general_eq13}
\left| \sum_{j,k,l} \varepsilon_{5,j,k,l} + \varepsilon_{6,j,k,l} \right| \le f(z) \E(Y_{1,1}^4) \sqrt{n} \Tr(MM^t)^{1/2} \, .
\end{equation}

Finally, combining relations from (\ref{freeness_general_eq1}) to (\ref{freeness_general_eq13}), we get
\begin{multline}
\left| \E G_{\mu_{XX^t}}(z) - \E G_{\mu_{\Xhat\Xhat^t}}(z) \right| \\
\le f(z) \left( \E|Y_{1,1}|^3 + \E(Y_{1,1}^4) \right) \left( \frac{1}{\sqrt{n}} + \frac{\Tr(MM^t)^{1/2}}{n} \right) \, .
\end{multline} \qed
\end{Proof}

We can now conclude the proof of the general case and obtain Theorem \ref{theorem1bis}. In fact, in Proposition \ref{freeness_general_step1}, we assumed that $\E(Y_{1,1}) = 0$, so we only have to remove this assumption.

\begin{Proof}
We recall that $\mathring{X} = X-\E(X)$ by definition. We also define $\gbar(z) = \E G_{\mu_{XX^t}}(z)$, $g_{\circ}(z) = \E G_{\mu_{\mathring{X}\mathring{X}^t}}(z)$, and $\widehat{g}(z) = \E G_{\mu_{\Xhat\Xhat^t}}(z)$. Using the notations in Lemma \ref{freeness_lemma}, we have
\begin{eqnarray*}
\lefteqn{   \left| \gbar(z) - (1-c\gbar(z)) G_{\mu_{MM^t}}(z(1-c\gbar(z))^2-(1-c)(1-c\gbar(z))) \right|   } \\
& \le & \left| \gbar(z) - g_{\circ}(z) \right| + \left| g_{\circ}(z) - \widehat{g}(z) \right| + \left| \widehat{g}(z) - \phi_{z,\mu_{MM^t}}(\widehat{g}(z),c) \right| \\
& & \quad  + \ \left| \phi_{z,\mu_{MM^t}}(\widehat{g}(z),c) - \phi_{z,\mu_{MM^t}}(g_{\circ}(z),c) \right| \\
& & \quad  + \ \left| \phi_{z,\mu_{MM^t}}(g_{\circ}(z),c) - \phi_{z,\mu_{MM^t}}(\gbar(z),c) \right| \\
& \le & (1+l_{s,t}) \left| \gbar(z) - g_{\circ}(z) \right| + (1+l_{s,t}) \left| g_{\circ}(z) - \widehat{g}(z) \right| \\
& & \quad  + \ \left| \widehat{g}(z) - \phi_{z,\mu_{MM^t}}(\widehat{g}(z),c) \right|
\end{eqnarray*}
for $s$ large enough and $t$ small enough by Lemma \ref{freeness_lemma}.\\
Since the matrix $X-\mathring{X} = \E(X)$ has rank at most 1, using the relations (\ref{dst}), (\ref{comparison_distances}), and (\ref{rank_covariance}), we have
\begin{displaymath}
\left| G_{\mu_{XX^t}}(z) - G_{\mu_{\mathring{X}\mathring{X}^t}}(z) \right| \le d_{s,t} \left( \mu_{XX^t} , \mu_{\mathring{X}\mathring{X}^t} \right) \le \dKS \left( \mu_{XX^t} , \mu_{\mathring{X}\mathring{X}^t} \right) \le \frac1n \, .
\end{displaymath}
Proposition \ref{freeness_Gaussian} (the Gaussian case) applied to $\Yhat$ and Proposition \ref{freeness_general_step1} (the centred case) applied to $\mathring{Y}$ permit to get finally
\begin{multline*}
\left| \gbar(z) - (1-c\gbar(z)) G_{\mu_{MM^t}}(z(1-c\gbar(z))^2-(1-c)(1-c\gbar(z))) \right| \\
\le \ (1+l_{s,t}).\frac1n + \ f(z) \left( \E |\mathring{Y_{1,1}}|^3 + \E(\mathring{Y_{1,1}}^4) \right) \left( \frac{1}{\sqrt{n}} + \frac{\Tr(MM^t)^{1/2}}{n} \right) \\
+ \ f(z) \left( |c_n-c| + \frac1n + \frac{\Tr(MM^t)^{1/2}}{n^{5/4}} \right) \, .
\end{multline*} \qed
\end{Proof}

\section{Large deviations} \label{section_deviations}

This section is devoted to the proof of Theorem \ref{theorem2}. In this section, $X \in \MM_{n,p}(\R)$ is a random matrix such that $c_n = \frac{n}{p} \to c \in (0,+\infty)$. Moreover, we assume that $\Var(X_{1,1})=1$ and that there exist $\alpha \in (0,2)$ and $a \in (0,+\infty]$ such that $X_{1,1} \in \Sz_{\alpha}(a)$ (see Definition \ref{S_alpha(a)}).\\

We define
\begin{displaymath}
\varepsilon(n) = \frac{1}{\log n}
\end{displaymath}
and we decompose the matrix $X$ as
\begin{equation} \label{decomposition_X}
\frac{X}{\sqrt{p}} = A+B+C+D \, ,
\end{equation}
where $A,B,C,D$ are the matrices defined by
\begin{displaymath}
\begin{array}{lr}
A_{j,k} = \dfrac{X_{j,k}}{\sqrt{p}} \1_{|X_{j,k}| < (\log n)^{2/\alpha}}
& B_{j,k} = \dfrac{X_{j,k}}{\sqrt{p}} \1_{(\log n)^{2/\alpha} \le |X_{j,k}| \le \varepsilon(n) \sqrt{p}} \\
C_{j,k} = \dfrac{X_{j,k}}{\sqrt{p}} \1_{\varepsilon(n) \sqrt{p} < |X_{j,k}| \le \varepsilon(n)^{-1} \sqrt{p}}
& D_{j,k} = \dfrac{X_{j,k}}{\sqrt{p}} \1_{\varepsilon(n)^{-1} \sqrt{p} < |X_{j,k}|} \, .
\end{array}
\end{displaymath}

Besides, we denote by $B_{s,t}(\mu,\delta)$ the ball with centre $\mu \in \PP(\R)$ and radius $\delta>0$ for the distance $d_{s,t}$.

\subsection{Exponential equivalences}

The goal of this subsection is to prove the following.

\begin{Proposition} \label{exponential_equivalence}
There exist $s,t>0$ such that the random distributions $\mu_{XX^t/p}$ and $\left( \sqrt{\mu_{CC^t}} \boxplus_c \sqrt{\mump} \right)^2$ are $d_{s,t}$-exponentially equivalent at scale $n^{1+\alpha/2}$ as $n \to +\infty$, i.e. for all $\delta>0$, we have
\begin{displaymath}
\lim_{n \to +\infty} \frac{1}{n^{1+\alpha/2}} \log \PR \left( d_{s,t} \left( \mu_{XX^t/p} , \left( \sqrt{\mu_{CC^t}} \boxplus_c \sqrt{\mump} \right)^2 \right) \ge \delta \right) = -\infty \, .
\end{displaymath}
\end{Proposition}

The strategy to prove Proposition \ref{exponential_equivalence} is similar to the one in \cite{BC}. First, we explain why the contributions of $B$ and $D$ for large deviations can be neglected (Lemmas \ref{exp_equiv_lemma1} and \ref{exp_equiv_lemma2}) and then, we show that the measures $\mu_{(A+C)(A+C)^t}$ and $\left( \sqrt{\mu_{CC^t}} \boxplus_c \sqrt{\mump} \right)^2$ are exponentially equivalent thanks to a conditioning and a coupling argument in which several tools are needed, such as the concentration property (\ref{concentration3}) and the asymptotic freeness result stated in Theorem \ref{theorem1}. From now on, we consider $s>2$ and $t>0$.

First, the contribution of $D$ is negligible.

\begin{Lemma} \label{exp_equiv_lemma1}
$\mu_{XX^t/p}$ et $\mu_{(A+B+C)(A+B+C)^t}$ are exponentially equivalent.
\end{Lemma}

The proof is very similar to what is done in \cite{BC}, the only difference being the use of (\ref{rank_covariance}) instead of (\ref{rank_inequality}). Therefore, it will not be repeated here.

The contribution of $B$ is also negligible.

\begin{Lemma} \label{exp_equiv_lemma2}
$\mu_{XX^t/p}$ et $\mu_{(A+C)(A+C)^t}$ are exponentially equivalent.
\end{Lemma}

\begin{Proof}
From Lemma \ref{exp_equiv_lemma1}, the triangle inequality, Lemma 1.2.15 in \cite{DZ}, and the inequality $d_{s,t} \le W_1 \le W_2$, it is sufficient to prove that for all $\delta>0$,
\begin{displaymath}
\lim_{n \to +\infty} \frac{1}{n^{1+\alpha/2}} \log \PR \left( W_2 \left( \mu_{(A+B+C)(A+B+C)^t} , \mu_{(A+C)(A+C)^t} \right) \ge \delta \right) = -\infty \ .
\end{displaymath}
From (\ref{HW_covariance}), which is the analogue of the Hoffman-Wielandt inequality (\ref{HW_inequality}) for covariance matrices, it is sufficient to check that for all $\delta>0$,
\begin{multline*}
\lim_{n \to +\infty} \frac{1}{n^{1+\alpha/2}} \log \PR \Bigg( \frac{2}{n^2} \Tr((A+B+C)(A+B+C)^t \\
+ \, (A+C)(A+C)^t) \Tr(BB^t) \ge \delta \Bigg) = -\infty \ .
\end{multline*}
Let $\delta>0$. We have
\begin{displaymath}
\Tr((A+C)(A+C)^t) \le \Tr((A+B+C)(A+B+C)^t) \le \Tr \left( \frac1p XX^t \right)
\end{displaymath}
using the decomposition (\ref{decomposition_X}). Thus,
\begin{equation} \label{exp_equiv_eq1}
\begin{array}{l}
\displaystyle{   \PR \left( \frac{2}{n^2} \Tr((A+B+C)(A+B+C)^t + (A+C)(A+C)^t) \Tr(BB^t) \ge \delta \right)   } \\
\displaystyle{   \ \le \ \PR \left( \frac{4}{n^2p} \Tr(XX^t) \Tr(BB^t) \ge \delta \right)   } \\
\displaystyle{   \ \le \ \PR \left( \frac{1}{np} \Tr(XX^t) \ge \E(X_{1,1}^2) + \delta \right)   } \\
\displaystyle{   \ \qquad + \, \PR \left( \frac4n \Tr(BB^t) \ge \frac{\delta}{\E(X_{1,1}^2) + \delta} \right)   } \, .
\end{array}
\end{equation}

On the one hand, since $\Tr(XX^t)$ is the sum of $np$ i.i.d. random variables, from Cram\'er's theorem in $\R$ (see \cite[Theorem 2.2.3]{DZ}), we have
\begin{multline*}
\lim_{n \to +\infty} \frac{1}{np} \log \PR \left( \frac{1}{np} \Tr(XX^t) \ge \E(X_{1,1}^2) + \delta \right) \\
= -\sup_{\theta \in \R} \left( \theta(\E(X_{1,1}^2) + \delta) - \log \E(e^{\theta X_{1,1}}) \right) <0
\end{multline*}
so, since $\alpha<2$,
\begin{equation} \label{exp_equiv_eq2}
\lim_{n \to +\infty} \frac{1}{n^{1+\alpha/2}} \log \PR \left( \frac{1}{np} \Tr(XX^t) \ge \E(X_{1,1}^2) + \delta \right) = -\infty \, .
\end{equation}

On the other hand, since $\frac{n}{p} \to c$, the same arguments as in \cite{BC} lead to
\begin{equation} \label{exp_equiv_eq3}
\lim_{n \to +\infty} \frac{1}{n^{1+\alpha/2}} \log \PR \left( \frac4n \Tr(BB^t) \ge \frac{\delta}{\E(X_{1,1}^2) + \delta} \right) = -\infty \, .
\end{equation}
Finally, combining (\ref{exp_equiv_eq1}), (\ref{exp_equiv_eq2}), (\ref{exp_equiv_eq3}), and Lemma 1.2.15 in \cite{DZ}, we get the exponential equivalence of $\mu_{XX^t/p}$ and $\mu_{(A+C)(A+C)^t}$. \qed
\end{Proof}

Before proving Proposition \ref{exponential_equivalence}, we need some additional properties.

\begin{Lemma} \label{exp_equiv_lemma3}
\begin{itemize}
\item[(i)] We have
\begin{displaymath}
\lim_{n \to +\infty} \frac{1}{n^{1+\alpha/2}} \log \PR \left( \frac1n \Tr(CC^t) > (\log n)^2 \right) = -\infty \, .
\end{displaymath}
\item[(ii)] Defining $I=\{ (j,k) \ | \ |X_{j,k}| \ge (\log n)^{2/\alpha} \}$, for all $\delta>0$, we have
\begin{displaymath}
\lim_{n \to +\infty} \frac{1}{n^{1+\alpha/2}} \log \PR \left( |I| \ge \delta n^{1+\alpha/2} \right) = -\infty \, .
\end{displaymath}
\item[(iii)] We denote by $P_n$ the distribution of $X_{1,1}$ conditionally to $\{ |X_{1,1}| < (\log n)^{2/\alpha} \}$. Let $Z_n$ be a random variable with distribution $P_n$. There exists $\zeta>0$ such that
\begin{displaymath}
\sup_{n \in \N} \max \left( \E(Z_n^2) , (\E(Z_n^2))^2 , \E(Z_n^4) \right) \le \zeta \, .
\end{displaymath}
Furthermore, the variance of $Z_n$, denoted by $\sigma_n^2$, tends to $\Var(X_{1,1})=1$ as $n \to +\infty$ and more precisely, there exists $\eta>0$ such that
\begin{displaymath}
|\sigma_n^2-1| \le \eta e^{-a(\log n)^2/4} \, .
\end{displaymath}
\end{itemize}
\end{Lemma}

\begin{Proof}
The proofs of (i) and (ii) exactly follow the proof of Lemma 2.4 in \cite{BC}. Therefore, we will only prove (iii).\\
Let $Z_n$ be a random variable with distribution $P_n$ defined as above. We have
\begin{displaymath}
\E(Z_n^2) = \E( X_{1,1}^2 \ | \ |X_{1,1}| < (\log n)^{2/\alpha} ) = \frac{\E \left( X_{1,1}^2 \1_{|X_{1,1}| < (\log n)^{2/\alpha}} \right)}{\PR(|X_{1,1}| < (\log n)^{2/\alpha})} \, .
\end{displaymath}
But thanks to hypothesis (\ref{hypo1}), $X_{1,1}^2$ is integrable, so by the dominated convergence theorem, $\E \left( X_{1,1}^2 \1_{|X_{1,1}| < (\log n)^{2/\alpha}} \right)$ tends to $\E(X_{1,1}^2)$ as $n \to +\infty$. Besides, $\PR(|X_{1,1}| < (\log n)^{2/\alpha})$ tends to 1, so $\E(Z_n^2)$ tends to $\E(X_{1,1}^2)$ as $n \to +\infty$.\\
The same arguments show that $\E(Z_n^4)$ tends to $\E(X_{1,1}^4)$ as $n \to +\infty$. We can deduce that there exists a real number $\zeta$ such that
\begin{displaymath}
\sup_{n \in \N} \max \left( \E(Z_n^2) , (\E(Z_n^2))^2 , \E(Z_n^4) \right) \le \zeta \, .
\end{displaymath}

Moreover, we have
\begin{eqnarray*}
\sigma_n^2 & = & \Var( X_{1,1} \ | \ |X_{1,1}| < (\log n)^{2/\alpha} ) \\
& = & \frac{\E \left( X_{1,1}^2 \1_{|X_{1,1}| < (\log n)^{2/\alpha}} \right)}{\PR(|X_{1,1}| < (\log n)^{2/\alpha})} - \left( \frac{\E \left( X_{1,1} \1_{|X_{1,1}| < (\log n)^{2/\alpha}} \right)}{\PR(|X_{1,1}| < (\log n)^{2/\alpha})} \right)^2 \, .
\end{eqnarray*}
Using similar arguments, we prove that $\sigma_n^2$ tends to $\Var(X_{1,1})=1$ as $n \to +\infty$. More precisely, we can write
\begin{eqnarray} \label{exp_equiv_eq3bis}
\sigma_n^2-1 & = & \frac{\E \left( X_{1,1}^2 \1_{|X_{1,1}| < (\log n)^{2/\alpha}} \right)}{\PR \left( |X_{1,1}| < (\log n)^{2/\alpha} \right)} - \left( \frac{\E \left( X_{1,1} \1_{|X_{1,1}| < (\log n)^{2/\alpha}} \right)}{\PR \left( |X_{1,1}| < (\log n)^{2/\alpha} \right)} \right)^2 \nonumber \\
& & \quad - \E(X_{1,1}^2) + (\E(X_{1,1}))^2 \nonumber \\
& = & \frac{\E \left( X_{1,1}^2 \1_{|X_{1,1}| < (\log n)^{2/\alpha}} \right) - \E(X_{1,1}^2) \PR \left( |X_{1,1}| < (\log n)^{2/\alpha} \right)}{\PR \left( |X_{1,1}| < (\log n)^{2/\alpha} \right)} \nonumber \\
& & \quad + \frac{\E(X_{1,1})^2 \PR \left( |X_{1,1}| < (\log n)^{2/\alpha} \right)^2 - \E \left( X_{1,1} \1_{|X_{1,1}| < (\log n)^{2/\alpha}} \right)^2}{\PR \left( |X_{1,1}| < (\log n)^{2/\alpha} \right)^2} \nonumber \\
& = & \frac{\E(X_{1,1}^2) \PR \left( |X_{1,1}| \ge (\log n)^{2/\alpha} \right) - \E \left( X_{1,1}^2 \1_{|X_{1,1}| \ge (\log n)^{2/\alpha}} \right)}{\PR \left( |X_{1,1}| < (\log n)^{2/\alpha} \right)} \nonumber \\
& & \quad + \frac{\E(X_{1,1})^2 \left( \PR \left( |X_{1,1}| < (\log n)^{2/\alpha} \right)^2 - 1 \right)}{\PR \left( |X_{1,1}| < (\log n)^{2/\alpha} \right)^2} \nonumber \\
& & \quad + \frac{2 \E(X_{1,1}) \E \left( X_{1,1} \1_{|X_{1,1}| \ge (\log n)^{2/\alpha}} \right)}{\PR \left( |X_{1,1}| < (\log n)^{2/\alpha} \right)^2} \nonumber \\
& & \quad - \frac{\E \left( X_{1,1} \1_{|X_{1,1}| \ge (\log n)^{2/\alpha}} \right)^2}{\PR \left( |X_{1,1}| < (\log n)^{2/\alpha} \right)^2}
\end{eqnarray}
where $\E \left( X_{1,1} \1_{|X_{1,1}| < (\log n)^{2/\alpha}} \right)^2 = \left( \E(X_{1,1}) - \E \left( X_{1,1} \1_{|X_{1,1}| \ge (\log n)^{2/\alpha}} \right) \right)^2$ was used to get the last equality.\\
From hypothesis (\ref{hypo1}), for $n$ large enough, we have
\begin{displaymath}
\PR(|X_{1,1}| \ge (\log n)^{2/\alpha}) \le e^{-\frac{a}{2} (\log n)^2} \, .
\end{displaymath}
Besides,
\begin{multline*}
\PR(|X_{1,1}| < (\log n)^{2/\alpha})^2 - 1 \\
= - \PR(|X_{1,1}| \ge (\log n)^{2/\alpha}) \left( \PR(|X_{1,1}| < (\log n)^{2/\alpha}) + 1 \right)
\end{multline*}
and by the Cauchy-Schwarz inequality, we have
\begin{displaymath}
\left| \E \left( X_{1,1} \1_{|X_{1,1}| \ge (\log n)^{2/\alpha}} \right) \right| \le \E(X_{1,1}^2)^{1/2} \PR(|X_{1,1}| \ge (\log n)^{2/\alpha})^{1/2}
\end{displaymath}
and
\begin{displaymath}
\left| \E \left( X_{1,1}^2 \1_{|X_{1,1}| \ge (\log n)^{2/\alpha}} \right) \right| \le \E(X_{1,1}^4)^{1/2} \PR(|X_{1,1}| \ge (\log n)^{2/\alpha})^{1/2} \, .
\end{displaymath}
Going back to (\ref{exp_equiv_eq3bis}), we have for $n$ large enough
\begin{multline*}
|\sigma_n^2-1| \le 2 \E(X_{1,1}^2) e^{-\frac{a}{2} (\log n)^2} + 2 \E(X_{1,1}^4)^{1/2} e^{-\frac{a}{4} (\log n)^2} \\
+ \ 4 \E(X_{1,1})^2 e^{-\frac{a}{2} (\log n)^2} + 4 |\E(X_{1,1})| \E(X_{1,1}^2)^{1/2} e^{-\frac{a}{4} (\log n)^2} \\
+ \ 2 \E(X_{1,1}^2) e^{-\frac{a}{2} (\log n)^2} \, .
\end{multline*}
Because the moments of $X_{1,1}$ are finite, we can deduce that there exists a real number $\eta$ such that
\begin{displaymath}
|\sigma_n^2-1| \le \eta e^{-a(\log n)^2/4} \, .
\end{displaymath} \qed
\end{Proof}

We can now prove Proposition \ref{exponential_equivalence}.

\begin{Proof}
The proof relies on a conditioning with respect to the entries of $X$ which are not in $A$ and on a coupling argument to remove the dependency between $A$ and $C$.

We use here the same notations as \cite{BC}. We denote by $\FF$ the $\sigma$-algebra
\begin{displaymath}
\FF = \sigma \left\{ X_{j,k} \1_{|X_{j,k}| \ge (\log n)^{2/\alpha}} \right\} \, ,
\end{displaymath}
$\PR_{\FF}$ and $\E_{\FF}$ the probability and the expectation conditionally to $\FF$, and we denote by $E$ and $F$ the events
\begin{displaymath}
E = \left\{ \frac1n \Tr(CC^t) \le (\log n)^2 \right\}
\end{displaymath}
and
\begin{displaymath}
F = \left\{ |I| < n^{1+\alpha/2} \right\} \, ,
\end{displaymath}
with $I=\{ (j,k) \ | \ |X_{j,k}| \ge (\log n)^{2/\alpha} \}$. Thus, the matrix $C$ is $\FF$-measurable and the events $E$ and $F$ belong to $\FF$. Moreover, from Lemma \ref{exp_equiv_lemma3} (i)-(ii), we have
\begin{equation} \label{exp_equiv_eq4}
\lim_{n \to +\infty} \frac{1}{n^{1+\alpha/2}} \log \PR(E^c) = -\infty \quad \mathrm{and} \quad \lim_{n \to +\infty} \frac{1}{n^{1+\alpha/2}} \log \PR(F^c) = -\infty \, .
\end{equation}

Besides, conditionally to $\FF$, $\sqrt{p}A$ is a random matrix with independent entries bounded by $(\log n)^{2/\alpha}$. From the concentration result (\ref{concentration3}) applied to $Y=\sqrt{p}A$, $M=C$, $\kappa=(\log n)^{2/\alpha}$, from the inequality $d_{s,t} \le W_1$, and using that $\alpha<2$, we get for all $\delta>0$ and $n$ large enough,
\begin{multline*}
\1_E \PR_{\FF} \left( d_{s,t} \left( \mu_{(A+C)(A+C)^t} , \E_{\FF} \mu_{(A+C)(A+C)^t} \right) \ge \delta \right) \\
\le \frac{\beta (\log n)^{2/\alpha}}{\delta^{3/2}} \exp \left( -\frac{n^2\delta^5}{\beta (\log n)^{8/\alpha}} \right)
\end{multline*}
hence
\begin{equation} \label{exp_equiv_eq5}
\lim_{n \to +\infty} \frac{1}{n^{1+\alpha/2}} \log \PR \left( E \, \cap \, \left\{ d_{s,t} \left( \mu_{(A+C)(A+C)^t} , \E_{\FF} \mu_{(A+C)(A+C)^t} \right) \ge \delta \right\} \right) = -\infty \, .
\end{equation}

We will now use a coupling argument. We consider an independent random matrix $Y$ whose entries are i.i.d. with distribution $P_n$ defined in Lemma \ref{exp_equiv_lemma3}, and we denote by $A'$ the matrix defined by
\begin{displaymath}
A'_{j,k} = \1_{(j,k) \notin I} A_{j,k} + \1_{(j,k) \in I} \frac{Y_{j,k}}{\sqrt{p}} \, .
\end{displaymath}
Consequently, $\sqrt{p}A'$ and $Y$ have the same distribution and are independent from $\FF$. In particular, we will use later that for all bounded continuous $f$, we have $\E_{\FF}(f(Y)) = \E(f(Y))$.

From the inequalities (\ref{HW_covariance}) and $d_{s,t} \le W_2$, we have
\begin{eqnarray*}
\lefteqn{   d_{s,t} (\mu_{(A+C)(A+C)^t} , \mu_{(A'+C)(A'+C)^t})^4   } \\
& \le & \frac{2}{n^2} \Tr((A-A')(A-A')^t) \Tr((A+C)(A+C)^t + (A'+C)(A'+C)^t) \\
& = & \frac{2}{n^2} \left( \sum_{j,k} \1_{(j,k) \in I} \frac{Y_{j,k}^2}{p} \right) \left( \sum_{j,k} (A_{j,k}+C_{j,k})^2 + (A'_{j,k}+C_{j,k})^2 \right) \\
& = & \frac{2}{n^2p} \left( \sum_{(j,k) \in I} Y_{j,k}^2 \right) \left( \sum_{j,k} A_{j,k}^2 + 2C_{j,k}^2 + (A'_{j,k})^2 \right) \\
& = & \frac{2}{n^2p} \left( \sum_{(j,k) \in I} Y_{j,k}^2 \right) \left( \sum_{j,k} A_{j,k}^2 + 2 \Tr(CC^t) + \sum_{j,k} A_{j,k}^2 + \sum_{(j,k) \in I} \frac{Y_{j,k}^2}{p} \right) \\
& \le & \frac{2}{n^2p} \left( 2 \left( np \frac{(\log n)^{4/\alpha}}{p} + \Tr(CC^t) \right) \left( \sum_{(j,k) \in I} Y_{j,k}^2 \right) + \frac1p \left( \sum_{(j,k) \in I} Y_{j,k}^2 \right)^2 \right) \, .
\end{eqnarray*}
With definition (\ref{dst}) of $d_{s,t}$ and conditional Jensen's inequality for the concave function $x \mapsto x^{1/4}$, we thus have
\begin{eqnarray*}
\lefteqn{  \1_E \1_F d_{s,t} \left( \E_{\FF} \mu_{(A+C)(A+C)^t} , \E_{\FF} \mu_{(A'+C)(A'+C)^t} \right)  } \\
& \le & \1_E \1_F \E_{\FF} d_{s,t} \left( \mu_{(A+C)(A+C)^t} , \mu_{(A'+C)(A'+C)^t} \right) \\
& \le & \left[ \frac{2 \1_E \1_F}{n^2p} \left( 2(n(\log n)^{4/\alpha} + \Tr(CC^t)) \E_{\FF} \left( \sum_{(j,k) \in I} Y_{j,k}^2 \right) \right. \right. \\
& & + \left. \left. \frac1p \E_{\FF} \left( \sum_{\substack{(j,k) \in I \\ (l,m) \in I}} Y_{j,k}^2 Y_{l,m}^2 \right) \right) \right]^{1/4}
\end{eqnarray*}
because $\1_E$, $\1_F$, and $\Tr(CC^t)$ are $\FF$-measurable. Since the events $\{(j,k) \in I\}$ are $\FF$-measurable and $Y$ is independent from $\FF$, we have
\begin{displaymath}
\E_{\FF} \left( \sum_{j,k} \1_{(j,k) \in I} Y_{j,k}^2 \right) = \sum_{j,k} \1_{(j,k) \in I} \E(Y_{j,k}^2) = |I|.\E(Y_{1,1}^2) \le \zeta |I|
\end{displaymath}
from Lemma \ref{exp_equiv_lemma3} (iii), and similarly
\begin{displaymath}
\E_{\FF} \left( \sum_{j,k,l,m} \1_{(j,k) \in I} \1_{(l,m) \in I} Y_{j,k}^2 Y_{l,m}^2 \right) \le \zeta |I|^2 \, .
\end{displaymath}
So we have
\begin{eqnarray*}
\lefteqn{  \1_E \1_F d_{s,t} \left( \E_{\FF} \mu_{(A+C)(A+C)^t} , \E_{\FF} \mu_{(A'+C)(A'+C)^t} \right)  } \\
& \le & \left[ \frac{2 \zeta \1_E \1_F}{n^2p} \left( 2(n(\log n)^{4/\alpha} + \Tr(CC^t)) |I| + \frac1p |I|^2 \right) \right]^{1/4} \\
& \le & \left[ \frac{2 \zeta c_n}{n^3} \left( 2 \left( n (\log n)^{4/\alpha} + n (\log n)^2 \right) n^{1+\alpha/2} + c_n n^{1+\alpha} \right) \right]^{1/4} \\
& \le & \left[ \frac{2 \zeta c_n}{n^3} . 3n (\log n)^{4/\alpha} n^{1+\alpha/2} \right]^{1/4} \\
& = & \left( 6 \zeta c_n \right)^{1/4} \frac{(\log n)^{1/\alpha}}{n^{1/4-\alpha/8}}
\end{eqnarray*}
for $n$ large enough (we used here the fact that $\frac{4}{\alpha}>2$). It follows that for all $\delta>0$,
\begin{multline} \label{exp_equiv_eq6}
\lim_{n \to +\infty} \frac{1}{n^{1+\alpha/2}} \log \PR \left( E \, \cap \, F \, \right. \\
\left. \cap \, \left\{ d_{s,t} \left( \E_{\FF} \mu_{(A+C)(A+C)^t} , \E_{\FF} \mu_{(A'+C)(A'+C)^t} \right) \ge \delta \right\} \right) = -\infty \, .
\end{multline}

In addition, we define $\sigma_n^2 = \Var(Y_{1,1})$ as in Lemma \ref{exp_equiv_lemma3} (iii). Since $C$ is $\FF$-measurable, $Y$ is independent from $\FF$, and $\frac{1}{\sigma_n^4} \E_{\FF} (Y_{1,1}^4) \le 2\zeta < +\infty$ for $n$ large enough, we can apply Theorem \ref{theorem1} to $Y/\sigma_n$ and $C$, conditionally to $\FF$. Therefore, for $n$ large enough, $s$ large enough, and $t$ small enough, we have
\begin{eqnarray*}
\lefteqn{  \1_E d_{s,t} \left( \E_{\FF} \mu_{\left( \frac{Y}{\sigma_n \sqrt{p}} + C \right) \left( \frac{Y}{\sigma_n \sqrt{p}} + C \right)^t} , \left( \sqrt{\mu_{CC^t}} \boxplus_c \sqrt{\mump} \right)^2 \right)  } \\
& \le & c_{s,t} \1_E \left( \frac{1}{\sigma_n^3} \E |\mathring{Y_{1,1}}|^3 + \frac{1}{\sigma_n^4} \E(\mathring{Y_{1,1}}^4) \right) \left( \frac{1}{\sqrt{n}} + \frac{\Tr(CC^t)^{1/2}}{n} \right) \\
& & \quad + \, c_{s,t} \1_E \left( |c_n-c| + \frac1n + \frac{\Tr(CC^t)^{1/2}}{n^{5/4}} \right) \\
& \le & c_{s,t} \left( \frac{8(\log n)^{6/\alpha}}{\sigma_n^3} + \frac{16(\log n)^{8/\alpha}}{\sigma_n^4} \right) \left( \frac{1}{\sqrt{n}} + \frac{\log n}{\sqrt{n}} \right) \\
& & \quad + \, c_{s,t} \left( |c_n-c| + \frac1n + \frac{\log n}{n^{3/4}} \right)
\end{eqnarray*}
using Jensen's inequality and the fact that for all $j \in \llbracket 1,n \rrbracket$ and $k \in \llbracket 1,p \rrbracket$, we have $|\mathring{Y_{j,k}}| = |Y_{j,k} - \E_{\FF}(Y_{j,k})| \le 2(\log n)^{2/\alpha}$. Therefore, for all $\delta>0$,
\begin{multline} \label{exp_equiv_eq7}
\lim_{n \to +\infty} \frac{1}{n^{1+\alpha/2}} \log \PR \left( E \, \cap \, \right. \\
\left. \left\{ d_{s,t} \left( \E_{\FF} \mu_{\left( \frac{Y}{\sigma_n \sqrt{p}} + C \right) \left( \frac{Y}{\sigma_n \sqrt{p}} + C \right)^t} , \left( \sqrt{\mu_{CC^t}} \boxplus_c \sqrt{\mump} \right)^2 \right) \ge \delta \right\} \right) = -\infty \, .
\end{multline}

To finish, from (\ref{HW_covariance}), we have
\begin{eqnarray*}
\lefteqn{  d_{s,t} \left( \mu_{\left( \frac{Y}{\sqrt{p}} + C \right) \left( \frac{Y}{\sqrt{p}} + C \right)^t} ,  \mu_{\left( \frac{Y}{\sigma_n \sqrt{p}} + C \right) \left( \frac{Y}{\sigma_n \sqrt{p}} + C \right)^t} \right)^4  } \\
& \le & \frac{2}{n^2} \Tr \left( \left( 1-\frac{1}{\sigma_n} \right)^2 \frac{YY^t}{p} \right) \\
& & \quad \times \Tr \left( \left( \frac{Y}{\sqrt{p}} + C \right) \left( \frac{Y}{\sqrt{p}} + C \right)^t + \left( \frac{Y}{\sigma_n \sqrt{p}} + C \right) \left( \frac{Y}{\sigma_n \sqrt{p}} + C \right)^t \right) \\
& \le & \frac{2}{n^2p} \left( 1-\frac{1}{\sigma_n} \right)^2 \Tr(YY^t) \left( \sum_{j,k} \left( \frac{Y_{j,k}}{\sqrt{p}} + C_{j,k} \right)^2 + \left( \frac{Y_{j,k}}{\sigma_n \sqrt{p}} + C_{j,k} \right)^2 \right) \\
& \le & \frac{2}{n^2p} \left( 1-\frac{1}{\sigma_n} \right)^2 \Tr(YY^t) \left( 4\Tr(CC^t) + \frac2p \left( 1+\frac{1}{\sigma_n^2} \right) \Tr(YY^t) \right)
\end{eqnarray*}
so, using conditional Jensen's inequality and doing as above, we get
\begin{eqnarray*}
\lefteqn{  \1_E \E_{\FF} d_{s,t} \left( \mu_{\left( \frac{Y}{\sqrt{p}} + C \right) \left( \frac{Y}{\sqrt{p}} + C \right)^t} , \mu_{\left( \frac{Y}{\sigma_n \sqrt{p}} + C \right) \left( \frac{Y}{\sigma_n \sqrt{p}} + C \right)^t} \right)  } \\
& \le & \left[ \frac{2}{n^2p} \left( 1-\frac{1}{\sigma_n} \right)^2 \left( 4n(\log n)^2 \E_{\FF} \Tr(YY^t) + \frac2p \left( 1+\frac{1}{\sigma_n^2} \right) \E_{\FF}(\Tr(YY^t)^2) \right) \right]^{1/4} \\
& \le & \left[ \frac{2}{n^2p} \left( \frac{\sigma_n^2-1}{\sigma_n(\sigma_n+1)} \right)^2 \left( 4n(\log n)^2.np\zeta + \frac2p \left( 1+\frac{1}{\sigma_n^2} \right).n^2p^2\zeta \right) \right]^{1/4} \\
& = & \left( \frac{\sigma_n^2-1}{\sigma_n(\sigma_n+1)} \right)^{1/2} \left( 8\zeta (\log n)^2 + 4\zeta \left( 1+\frac{1}{\sigma_n^2} \right) \right)^{1/4} \, .
\end{eqnarray*}
By Lemma \ref{exp_equiv_lemma3} (iii), we deduce from it that for all $\delta>0$,
\begin{multline} \label{exp_equiv_eq8}
\lim_{n \to +\infty} \frac{1}{n^{1+\alpha/2}} \log \PR \left( E \, \cap \, \right. \\
\left. \left\{ d_{s,t} \left( \E_{\FF} \mu_{\left( \frac{Y}{\sqrt{p}} + C \right) \left( \frac{Y}{\sqrt{p}} + C \right)^t} , \E_{\FF} \mu_{\left( \frac{Y}{\sigma_n \sqrt{p}} + C \right) \left( \frac{Y}{\sigma_n \sqrt{p}} + C \right)^t} \right) \ge \delta \right\} \right) = -\infty \, .
\end{multline}

To conclude, combining equalities from (\ref{exp_equiv_eq4}) to (\ref{exp_equiv_eq8}), Lemma \ref{exp_equiv_lemma2}, and Lemma 1.2.15 in \cite{DZ}, for $s$ large enough and $t$ small enough, we have for all $\delta>0$,
\begin{displaymath}
\lim_{n \to +\infty} \frac{1}{n^{1+\alpha/2}} \log \PR \left( d_{s,t} \left( \mu_{XX^t/p} , \left( \sqrt{\mu_{CC^t}} \boxplus_c \sqrt{\mump} \right)^2 \right) \ge \delta \right) = -\infty \, .
\end{displaymath}
\qed
\end{Proof}

\subsection{Large deviations for $\mu_{C'}$}

In the previous subsection, we proved that $\mu_{XX^t/p}$ and $\left( \sqrt{\mu_{CC^t}} \boxplus_c \sqrt{\mump} \right)^2$ are exponentially equivalent. Consequently, to obtain the large deviations of $\mu_{XX^t/p}$ (Theorem \ref{theorem2}), it is sufficient to study the large deviations of $\mu_{CC^t}$ and to apply the contraction principle (see \cite[Theorem 4.2.1]{DZ}). For this, in this subsection, we will study the large deviations of
\begin{displaymath}
C' = \left( \begin{array}{c|c} 0 & C \\ \hline C^t & 0 \end{array} \right)
\end{displaymath}
and prove the following, from which we will deduce the large deviations of $\mu_{CC^t}$ thanks to identity (\ref{muM'2}) and conclude in the next subsection.

\begin{Proposition} \label{LDP_C'}
The measure $\mu_{C'}$ satisfies the LDP with speed $n^{1+\alpha/2}$ in $\PP(\R)$, for weak topology and good rate function $\Phi'$ defined by
\begin{equation} \label{Phi'}
\Phi'(\mu) = \left\{ \begin{array}{ll} \frac{a}{2} \frac{c+1}{c^{1+\alpha/2}} m_{\alpha}(\mu) & \textrm{if $\mu$ is symmetric and } \mu(\{0\}) \ge \frac{|1-c|}{1+c} \\ +\infty & \textrm{otherwise} \end{array} \right. \, ,
\end{equation}
where $m_{\alpha}(\mu) = \int_{\R} |x|^{\alpha} \, d\mu(x)$ denotes the $\alpha$-th moment of $\mu$.
\end{Proposition}

Note that $\Phi'$ is a good rate function because it is well known that for all $m \ge 0$ and $p>0$, the set
\begin{equation} \label{Kpm}
K_{p,m} = \left\{ \mu \in \PP(\R) \ \left| \ \int_{\R} |x|^p \, d\mu(x) \le m \right. \right\}
\end{equation}
is compact for the weak topology. Moreover, the domain of $\Phi'$ can be explained thanks to Lemma \ref{LDP_C'_lemma1} (i).

\begin{Lemma} \label{LDP_C'_lemma1}
Let $M \in \MM_{n,p}(\R)$ and
\begin{displaymath}
M' = \left( \begin{array}{c|c} 0 & M \\ \hline M^t & 0 \end{array} \right) \, .
\end{displaymath}
\begin{itemize}
\item[(i)] The distribution $\mu_{M'}$ is symmetric and $\mu_{M'}(\{0\}) \ge \frac{|1-c_n|}{1+c_n}$.
\item[(ii)] We have
\begin{equation} \label{muM'2}
\mu_{M'}^2 = \frac{2c_n}{c_n+1} \mu_{MM^t} + \frac{1-c_n}{1+c_n} \delta_0 \, .
\end{equation}
\item[(iii)] If $M$ is diagonal, in the sense that only the entries $M_{j,j}$, $1 \le j \le n \wedge p$, can be non-zero, then
\begin{equation} \label{muM'}
\mu_{M'} = \frac{1}{n+p} \sum_{j=1}^{n \wedge p} \left( \delta_{M_{j,j}} + \delta_{-M_{j,j}} \right) + \frac{|1-c_n|}{1+c_n} \delta_0 \, .
\end{equation}
\end{itemize}
\end{Lemma}

The proof of this lemma does not present any difficulty and is left to the reader. We also need a second lemma, which consists in two estimates for the distribution of $X_{1,1}$. These estimates come from the particular form of this distribution, see hypotheses (\ref{hypo1}) and (\ref{hypo2bis}).

\begin{Lemma}
\begin{itemize}
\item[(i)] There exists a sequence $(\eta_n)_{n \in \N}$ converging to 0 such that for all $x \ge \varepsilon(n)$, we have
\begin{equation} \label{LDP_C'_lemma2}
\PR(|X_{1,1}| \ge x\sqrt{p}) \le e^{-(a-\eta_n) x^{\alpha} p^{\alpha/2}} \, .
\end{equation}
\item[(ii)] We denote by $S_a$ the support of the distribution $\vartheta_a$ defined by (\ref{hypo2bis}). There exists a sequence $(a_n)_{n \in \N}$ converging to $a$ such that for all $x \in \R$ satisfying $|x| \ge \varepsilon(n)$ and $\sgn(x) \in S_a$, for all $\gamma>0$, and for all $n$ large enough, we have
\begin{equation} \label{LDP_C'_lemma3}
\PR \left( \frac{X_{1,1}}{\sqrt{p}} \in (x-\gamma,x+\gamma) \right) \ge e^{-a_n |x|^{\alpha} p^{\alpha/2}} \, .
\end{equation}
\end{itemize}
\end{Lemma}

The computations to get these inequalities are explained in \cite[p. 26]{BC} and are left to the reader.

We will now prove Proposition \ref{LDP_C'}. Let us mention that Schatten's inequality (\ref{Schatten}) will be crucial in the proof since it will allow to link the $\alpha$-th moment of the spectral measure $\mu_{C'}$ and the entries of $C'$.

\begin{Proof}
Since the set of symmetric probability measures on $\R$ is closed for the weak topology, it is enough to prove the LDP on this set, see \cite[Lemma 4.1.5]{DZ}.\\

\emph{Upper bound.} Let $\mu$ be a symmetric probability measure on $\R$. Since the function $m_{\alpha}$ is lower semi-continuous, there exists a continuous function $h$ such that $h(0)=0$ and
\begin{displaymath}
\PR(\mu_{C'} \in B_{s,t}(\mu,\delta)) \le \PR(m_{\alpha}(\mu_{C'}) \ge m_{\alpha}(\mu) - h(\delta))
\end{displaymath}
for all $\delta$ small enough. Moreover, by Schatten's inequality (\ref{Schatten}) and the fact that $\sum_{i=1}^k a_i^r \le \left( \sum_{i=1}^k a_i \right)^r$ for all $r \ge 1$, $a_1,\ldots,a_k \ge 0$, we have
\begin{displaymath}
m_{\alpha}(\mu_{C'}) \le \frac{1}{n+p} \sum_{j=1}^{n+p} \left( \sum_{k=1}^{n+p} (C'_{j,k})^2 \right)^{\alpha/2} \le \frac{1}{n+p} \sum_{j=1}^{n+p} \sum_{k=1}^{n+p} |C'_{j,k}|^{\alpha} \, .
\end{displaymath}
Consequently,
\begin{eqnarray*}
\lefteqn{   \PR(\mu_{C'} \in B_{s,t}(\mu,\delta))   } \\
& \le & \PR \left( \frac{1}{n+p} \sum_{j,k} |C'_{j,k}|^{\alpha} \ge m_{\alpha}(\mu) - h(\delta) \right) \\
& = & \PR \left( \frac{2}{(c_n+1)p^{1+\alpha/2}} \sum_{j,k} |X_{j,k}|^{\alpha} \1_{\varepsilon(n) \sqrt{p} < |X_{j,k}| \le \varepsilon(n)^{-1} \sqrt{p}} \ge m_{\alpha}(\mu) - h(\delta) \right) \\
& \le & e^{-\frac{a_1}{2} (c_n+1) p^{1+\alpha/2} (m_{\alpha}(\mu) - h(\delta))} \left( \E \left( e^{a_1 |X_{1,1}|^{\alpha} \1_{\varepsilon(n) \sqrt{p} < |X_{j,k}| \le \varepsilon(n)^{-1} \sqrt{p}}} \right) \right)^{np}
\end{eqnarray*}
for all $a_1 \in (0,a)$ by Chernoff's inequality. Besides, from hypothesis (\ref{hypo1}), there exists $a_2 \in (a_1,a)$ such that for all $x$ large enough, $\PR(|X_{1,1}| \ge x) \le \exp(-a_2 x^{\alpha})$. Let us also recall the following integration by parts formula: for all $\nu \in \PP(\R)$ and $f \in \CC^1(\R,\R)$,
\begin{displaymath}
\int_a^b f(x) \, d\nu(x) = f(a) \nu([a,+\infty)) - f(b) \nu([b,+\infty)) + \int_a^b f'(x) \nu([x,+\infty)) \, dx \, .
\end{displaymath}

Denoting by $\PR_{|X_{1,1}|}$ the law of $|X_{1,1}|$, we thus have
\begin{eqnarray*}
\lefteqn{   \E \exp \left( a_1 |X_{1,1}|^{\alpha} \1_{\varepsilon(n) \sqrt{p} < |X_{1,1}| \le \varepsilon(n)^{-1} \sqrt{p}} \right)   } \\
& \le & 1 + \int_{\varepsilon(n) \sqrt{p}}^{\varepsilon(n)^{-1} \sqrt{p}} e^{a_1 x^{\alpha}} \, d\PR_{|X_{1,1}|}(x) \\
& \le & 1 + e^{a_1 \varepsilon(n)^{\alpha} p^{\alpha/2} - a_2 \varepsilon(n)^{\alpha} p^{\alpha/2}} + \int_{\varepsilon(n) \sqrt{p}}^{\varepsilon(n)^{-1} \sqrt{p}} a_1 \alpha x^{\alpha-1} e^{a_1 x^{\alpha} - a_2 x^{\alpha}} \, dx \\
& \le & 1 + e^{-(a_2-a_1) \varepsilon(n)^{\alpha} p^{\alpha/2}} - \frac{a_1}{a_1-a_2} e^{-(a_2-a_1) \varepsilon(n)^{\alpha} p^{\alpha/2}} \\
& = & 1 + \frac{a_2}{a_2-a_1} e^{-(a_2-a_1) \varepsilon(n)^{\alpha} p^{\alpha/2}} \\
& \le & \exp \left( \frac{a_2}{a_2-a_1} e^{-(a_2-a_1) \varepsilon(n)^{\alpha} p^{\alpha/2}} \right)
\end{eqnarray*}
hence
\begin{multline*}
\PR(\mu_{C'} \in B_{s,t}(\mu,\delta)) \\
\le \exp \left( -\frac{a_1}{2} (c_n+1) p^{1+\alpha/2} (m_{\alpha}(\mu)-h(\delta)) + \frac{a_2}{a_2-a_1} np e^{-(a_2-a_1) \varepsilon(n)^{\alpha} p^{\alpha/2}} \right) \, .
\end{multline*}
So, for all $\delta$ small enough and all $a_1 \in (0,a)$,
\begin{displaymath}
\limsup_{n \to +\infty} \frac{1}{n^{1+\alpha/2}} \log \PR(\mu_{C'} \in B_{s,t}(\mu,\delta)) \le -\frac{a_1}{2} \frac{c+1}{c^{1+\alpha/2}} (m_{\alpha}(\mu) - h(\delta))
\end{displaymath}
and finally
\begin{equation} \label{LDP_C'_eq1}
\limsup_{\delta \to 0} \limsup_{n \to +\infty} \frac{1}{n^{1+\alpha/2}} \log \PR(\mu_{C'} \in B_{s,t}(\mu,\delta)) \le -\frac{a}{2} \frac{c+1}{c^{1+\alpha/2}} m_{\alpha}(\mu) \, .
\end{equation}

In the case of a $\mu$ satisfying $\mu(\{0\}) < \frac{|1-c|}{1+c}$, we have a better result. Indeed, inspired by \cite{HP}, we can observe that for all $\varepsilon$ small enough, there exists $R>0$ such that
\begin{displaymath}
\mu([-R,R]) < \frac{|1-c|}{1+c} - \varepsilon \, .
\end{displaymath}
Therefore,
\begin{displaymath}
\left\{ \mu' \in \PP(\R) \ | \ \mu'([-R,R]) < \frac{|1-c|}{1+c} - \varepsilon \right\}
\end{displaymath}
is a neighbourhood of $\mu$ in which, for $n$ large enough, almost surely, $\mu_{C'}$ is not. So we have
\begin{equation} \label{LDP_C'_eq2}
\lim_{\delta \to 0} \lim_{n \to +\infty} \frac{1}{n^{1+\alpha/2}} \log \PR(\mu_{C'} \in B_{s,t}(\mu,\delta)) = -\infty \, .
\end{equation}

We have obtained the upper bound of the LDP.\\

\emph{Lower bound.} Let $\mu \in \PP(\R)$ be a symmetric measure such that $\mu(\{0\}) \ge \frac{|1-c|}{1+c}$. There exists $\tilde{\mu} \in \PP(\R_+)$ such that
\begin{displaymath}
\mu = \frac{|1-c|}{1+c} \delta_0 + \frac{1 \wedge c}{1+c} (\tilde{\mu} + (-\Id)_{\sharp}\tilde{\mu}) \, ,
\end{displaymath}
where $(-\Id)_{\sharp}\tilde{\mu}$ denotes the push-forward of $\tilde{\mu}$ by $-\Id$.\\
We denote by $x_1,\ldots,x_{n \wedge p}$ the quantiles of $\tilde{\mu}$ of orders $\frac{1}{1+n \wedge p}, \ldots, \frac{n \wedge p}{1+n \wedge p}$, we also define $n_0 = \min \{ j \in \llbracket 1,n \wedge p \rrbracket \ | \ x_j \ge \varepsilon(n) \}$, and
\begin{displaymath}
M' = \left( \begin{array}{c|c} 0 & M \\ \hline M^t & 0 \end{array} \right)
\end{displaymath}
with $M \in \MM_{n,p}(\R)$ defined by $M_{j,j} = x_j$ for all $j \in \llbracket n_0, n \wedge p \rrbracket$ and $M_{j,k}=0$ otherwise.\\
From (\ref{muM'}), we have
\begin{displaymath}
\mu_{M'} = \frac{1}{n+p} \sum_{j=1}^{n \wedge p} \left( \delta_{M_{j,j}} + \delta_{-M_{j,j}} \right) + \frac{|1-c_n|}{1+c_n} \delta_0 \, .
\end{displaymath}
Besides,
\begin{equation} \label{LDP_C'_eq3}
m_{\alpha}(\mu) = \frac{2(1 \wedge c)}{1+c} \int_0^{+\infty} |x|^{\alpha} \, d\tilde{\mu}(x) \ge \frac{2(1 \wedge c)}{(1+c)(1+n \wedge p)} \sum_{j=1}^{n \wedge p} |M_{j,j}|^{\alpha} \, .
\end{equation}
Let us also remark that by construction, $d_{s,t}(\mu,\mu_{M'})$ tends to 0 as $n \to +\infty$.\\

Let $\delta>0$. For $n$ large enough, we thus have
\begin{equation} \label{LDP_C'_eq4}
d_{s,t}(\mu,\mu_{M'}) < \frac{\delta}{2} \, .
\end{equation}
Using (\ref{LDP_C'_eq4}), the Hoffman-Wielandt inequality (\ref{HW_inequality}), the independence of the $X_{j,k}$'s, the inequalities (\ref{LDP_C'_lemma2}) and (\ref{LDP_C'_lemma3}), the fact that $1 \le n_0 \le n \wedge p$, and (\ref{LDP_C'_eq3}), we get for $n$ large enough,
\begin{eqnarray*}
\lefteqn{   \PR(\mu_{C'} \in B_{s,t}(\mu,\delta))   } \\
& \ge & \PR \left( \mu_{C'} \in B_{s,t} \left( \mu_{M'},\frac{\delta}{2} \right) \right) \\
& \ge & \PR \left( \frac{1}{n+p} \Tr((C'-M')^2) \le \frac{\delta^2}{4} \right) \\
& = & \PR \left( \sum_{j,k} (C_{j,k}-M_{j,k})^2 \le \frac{\delta^2(n+p)}{8} \right) \\
& \ge & \PR \Bigg( \forall j \in \llbracket n_0,n \wedge p \rrbracket, \ (C_{j,j}-M_{j,j})^2 \le \frac{\delta^2(n+p)}{8(n \wedge p)} \\
& & \quad \quad \ \cap \ \ \forall (j,k) \ \textrm{different,} \ C_{j,k}=0 \Bigg) \\
& \ge & \PR \Bigg( \forall j \in \llbracket n_0,n \wedge p \rrbracket, \ \frac{X_{j,j}}{\sqrt{p}} \in \left( M_{j,j} - \sqrt{\frac{\delta^2(n+p)}{8(n \wedge p)}} , M_{j,j} + \sqrt{\frac{\delta^2(n+p)}{8(n \wedge p)}} \right) \\
& & \quad \quad \ \cap \ \ \forall (j,k) \ \textrm{different,} \ |X_{j,k}| < \varepsilon(n) \sqrt{p} \Bigg) \\
& \ge & \prod_{j=n_0}^{n \wedge p} e^{-a_n |M_{j,j}|^{\alpha} p^{\alpha/2}} \left( 1-e^{-(a-\eta_n) \varepsilon(n)^{\alpha} p^{\alpha/2}} \right)^{np - (n \wedge p - n_0 + 1)} \\
& \ge & \frac12 \exp \left( -a_n \sum_{j=n_0}^{n \wedge p} |M_{j,j}|^{\alpha} p^{\alpha/2} \right) \\
& \ge & \frac12 \exp \left( -\frac{a_n}{2} p^{1+\alpha/2} \frac{(1+c) \left( \frac1p + c \wedge 1 \right)}{1 \wedge c} m_{\alpha}(\mu) \right)
\end{eqnarray*}
Note that we can apply (\ref{LDP_C'_lemma2}) and (\ref{LDP_C'_lemma3}), even if it means to swap $\tilde{\mu}$ and $(-\Id)_{\sharp} \tilde{\mu}$ in order to apply (\ref{LDP_C'_lemma3}).
We finally get
\begin{equation} \label{LDP_C'_eq5}
\liminf_{n \to +\infty} \frac{1}{n^{1+\alpha/2}} \log \PR(\mu_{C'} \in B_{s,t}(\mu,\delta)) \ge -\frac{a}{2} \frac{c+1}{c^{1+\alpha/2}} m_{\alpha}(\mu)
\end{equation}
for all $\delta>0$.

This is the lower bound of the LDP.\\

\emph{Exponential tightness.} Let $A>0$ and $m = \frac{2Ac^{1+\alpha/2}}{a(1+c)}$. We recall that the set $K_{\alpha,m}$ defined by (\ref{Kpm}) is compact. Moreover, using the computations done for the upper bound, we have
\begin{multline*}
\PR(\mu_{C'} \notin K_{\alpha,m}) = \PR(m_{\alpha}(\mu_{C'}) > m) \\
\le \exp \left( -\frac{a_1}{2} (c_n+1) p^{1+\alpha/2} m + \frac{a_2}{a_2-a_1} np e^{-(a_2-a_1) \varepsilon(n)^{\alpha} p^{\alpha/2}} \right)
\end{multline*}
for all $a_1 \in (0,a)$ and some $a_2 \in (a_1,a)$. It follows that
\begin{equation} \label{LDP_C'_eq6}
\limsup_{n \to +\infty} \frac{1}{n^{1+\alpha/2}} \log \PR(\mu_{C'} \notin K_{\alpha,m}) \le -\frac{a}{2} \frac{c+1}{c^{1+\alpha/2}} m = -A \, .
\end{equation}

The combination of (\ref{LDP_C'_eq1}), (\ref{LDP_C'_eq2}), (\ref{LDP_C'_eq5}), and (\ref{LDP_C'_eq6}) is the desired LDP. \qed
\end{Proof}

\subsection{Conclusion}

To conclude this section, we show how to deduce the LDP for $\mu_{XX^t/p}$ (Theorem \ref{theorem2}) from the LDP for $\mu_{C'}$ (Proposition \ref{LDP_C'}).

\begin{Proposition} \label{LDP_CC^t}
The measure $\mu_{CC^t}$ satisfies the LDP with speed $n^{1+\alpha/2}$ in $\PP(\R_+)$, for the weak topology and the good rate function $\Psi'$ defined by
\begin{equation} \label{Psi}
\Psi'(\nu) = \left\{ \begin{array}{ll} \frac{a}{c^{\alpha/2}} m_{\alpha/2}(\nu) & \textrm{if } \nu(\{0\}) \ge \max \left( 0,1-\frac1c \right) \\ +\infty & \textrm{otherwise} \end{array} \right. \, .
\end{equation}
\end{Proposition}

\begin{Proof}
We define
\begin{displaymath}
T_n : \mu \mapsto \frac12 \left( 1+\frac{1}{c_n} \right) \mu^2 + \frac12 \left( 1-\frac{1}{c_n} \right) \delta_0
\end{displaymath}
and
\begin{displaymath}
T : \mu \mapsto \frac12 \left( 1+\frac1c \right) \mu^2 + \frac12 \left( 1-\frac1c \right) \delta_0 \, ,
\end{displaymath}
so that $\mu_{CC^t} = T_n(\mu_{C'})$, see (\ref{muM'2}).\\

Besides, we have
\begin{equation} \label{LDP_CC^t_eq1}
\lim_{n \to +\infty} d_{s,t}(\mu_{CC^t}, T(\mu_{C'})) = 0 \, .
\end{equation}
Indeed, let $n \in \N$ and $z$ in the domain $V_{s,t}$ defined by (\ref{Vst}). We have
\begin{eqnarray*}
\lefteqn{   \left| G_{\mu_{CC^t}}(z) - G_{T(\mu_{C'})}(z) \right|   } \\
& = & \left| \int_{\R} \frac{1}{z-x} \, d(T_n(\mu_{C'}))(x) - \int_{\R} \frac{1}{z-x} \, d(T(\mu_{C'}))(x) \right| \\
& = & \left| \frac12 \left( \frac{1}{c_n} - \frac1c \right) \int_{\R} \frac{1}{z-x} \, d(\mu_{C'}^2)(x) + \frac12 \left( \frac1c - \frac{1}{c_n} \right) \frac1z \right| \\
& \le & \frac{1}{|\Imag z|} \left| \frac{1}{c_n} - \frac1c \right| \, ,
\end{eqnarray*}
so, taking the upper bound on $z \in V_{s,t}$ and the limit as $n \to +\infty$, we get (\ref{LDP_CC^t_eq1}).\\

The contraction principle applied to the function $T$, see \cite[Theorem 4.2.1]{DZ}, will allow us to conclude. Indeed, $T$ takes its values in $\PP(\R_+)$ and is continuous for the weak topology. This strategy will make appear the good rate function $\Psi'$ defined for all $\nu \in \PP(\R_+)$ by
\begin{displaymath}
\Psi'(\nu) = \inf \{ \Phi'(\mu), \ \ \mu \in \PP(\R) \ \textrm{s.t.} \ \nu = T(\mu) \} \, .
\end{displaymath}
For all $\mu \in \PP(\R)$, we have $(T(\mu))(\{0\}) = \frac12 \left( 1+\frac1c \right) \mu(\{0\}) + \frac12 \left( 1-\frac1c \right)$, hence
\begin{displaymath}
\mu(\{0\}) \ge \frac{|1-c|}{1+c} \quad \Leftrightarrow \quad (T(\mu))(\{0\}) \ge \frac12 \frac{|1-c|+(c-1)}{c} = \max \left( 0,1-\frac1c \right) \, .
\end{displaymath}
Therefore, for all $\nu \in \PP(\R_+)$ such that $\nu(\{0\}) \ge \max \left( 0,1-\frac1c \right)$, there exists a symmetric $\mu \in \PP(\R)$ satisfying $\mu(\{0\}) \ge \frac{|1-c|}{1+c}$ and $\nu = T(\mu)$. We have in this case
\begin{displaymath}
\Phi'(\mu) = \frac{a}{2} \frac{c+1}{c^{1+\alpha/2}} \int_{\R} |x|^{\alpha} \, d\mu(x) = \frac{a}{c^{\alpha/2}} \int_{\R} |x|^{\alpha/2} \, d(T(\mu))(x) = \frac{a}{c^{\alpha/2}} m_{\alpha/2}(\nu)
\end{displaymath}
hence $\Psi'(\nu) = \frac{a}{c^{\alpha/2}} m_{\alpha/2}(\nu)$. In the case of $\nu(\{0\}) < \max \left( 0,1-\frac1c \right)$, we can not find a symmetric $\mu \in \PP(\R)$ satisfying $\mu(\{0\}) \ge \frac{|1-c|}{1+c}$ and $\nu = T(\mu)$, so $\Psi'(\nu) = +\infty$. Thus, we have computed $\Psi'(\nu)$ for every $\nu \in \PP(\R_+)$.\\

\emph{Lower bound.} Let $\mu \in \PP(\R_+)$ and $\delta>0$. From (\ref{LDP_CC^t_eq1}), for $n$ large enough, we have $d_{s,t}(\mu_{CC^t},T(\mu_{C'})) \le \frac{\delta}{2}$, hence
\begin{displaymath}
\PR(\mu_{CC^t} \in B_{s,t}(\mu,\delta)) \ge \PR \left( T(\mu_{C'}) \in B_{s,t} \left( \mu,\frac{\delta}{2} \right) \right) \, .
\end{displaymath}
By Proposition \ref{LDP_C'} and the contraction principle, we thus have
\begin{equation} \label{LDP_CC^t_eq2}
\liminf_{n \to +\infty} \frac{1}{n^{1+\alpha/2}} \log \PR(\mu_{CC^t} \in B_{s,t}(\mu,\delta)) \ge - \inf_{\nu \in B_{s,t}(\mu,\delta/2)} \Psi'(\nu) \ge -\Psi'(\mu) \, .
\end{equation}\\

\emph{Upper bound.} Let $F$ be a closed subset of $\PP(\R_+)$ and $\delta>0$. From (\ref{LDP_CC^t_eq1}), for $n$ large enough, we have $d_{s,t}(\mu_{CC^t},T(\mu_{C'})) \le \delta$, so
\begin{displaymath}
\PR(\mu_{CC^t} \in F) \le \PR(T(\mu_{C'}) \in F^{\delta}) \, ,
\end{displaymath}
where $F^{\delta}$ denotes the $\delta$-neighbourhood of $F$ for the distance $d_{s,t}$, namely
\begin{displaymath}
F^{\delta} = \{ \nu \in \PP(\R_+) \ | \ \exists \mu \in F, \, d_{s,t}(\mu,\nu) \le \delta \} \, .
\end{displaymath}
Applying the contraction principle again, we thus have
\begin{displaymath}
\limsup_{n \to +\infty} \frac{1}{n^{1+\alpha/2}} \log \PR(\mu_{CC^t} \in F) \le -\inf_{\nu \in F^{\delta}} \Psi'(\nu) \, .
\end{displaymath}
This is true for all $\delta>0$ so, taking the limit as $\delta \to 0$, we get (see \cite[Lemma 4.1.6(a)]{DZ})
\begin{equation} \label{LDP_CC^t_eq3}
\limsup_{n \to +\infty} \frac{1}{n^{1+\alpha/2}} \log \PR(\mu_{CC^t} \in F) \le -\inf_{\nu \in F} \Psi'(\nu) \, .
\end{equation}

Combining (\ref{LDP_CC^t_eq2}) and (\ref{LDP_CC^t_eq3}), we can conclude that $\mu_{CC^t}$ satisfies the announced LDP. \qed
\end{Proof}

Because rectangular free convolution is continuous for weak topology, see \cite[Theorem 3.12]{BG}, the function $\mu \mapsto \left( \sqrt{\mu} \boxplus_c \sqrt{\mump} \right)^2$ is so, therefore, by Proposition \ref{LDP_CC^t} and the contraction principle, $\left( \sqrt{\mu_{CC^t}} \boxplus_c \sqrt{\mump} \right)^2$ satisfies the LDP with speed $n^{1+\alpha/2}$ on $\PP(\R_+)$ governed by the good rate function
\begin{displaymath}
J'(\mu) = \left\{ \begin{array}{ll} \Psi'(\nu) & \textrm{ if there exists } \nu \in \PP(\R_+) \textrm{ such that } \mu = \left( \sqrt{\nu} \boxplus_c \sqrt{\mump} \right)^2 \\ +\infty & \textrm{ otherwise} \end{array} \right. \, .
\end{displaymath}
Thanks to the exponential equivalence between $\mu_{XX^t/p}$ and $\left( \sqrt{\mu_{CC^t}} \boxplus_c \sqrt{\mump} \right)^2$ obtained in Proposition \ref{exponential_equivalence}, we can conclude that $\mu_{XX^t/p}$ satisfies the same LDP, see \cite[Theorem 4.2.13]{DZ}, which ends the proof of Theorem \ref{theorem2}.

\appendix

\section{Concentration bounds for the information-plus-noise model} \label{appendix_concentration}

\subsection{Concentration for some functions of the resolvent}

In Section \ref{section_freeness}, in order to prove Lemma \ref{freeness_Gaussian_lemma2}, we needed the following concentration estimates.

\begin{Proposition}[adaptation from {\cite[Lemma 8]{VLM}}] \label{concentration1-2}
Let $Y \in \MM_{n,p}(\R)$ be a random matrix with i.i.d. entries, let $M \in \MM_{n,p}(\R)$ be a deterministic matrix and let $z \in \C \setminus \R$.\\
We define $X = \frac{Y}{\sqrt{p}}+M$, $S=(zI_n-XX^t)^{-1}$ the resolvent of $XX^t$, $c_n = \frac{n}{p}$, and $\sigma^2 = \Var(Y_{1,1})$.\\
We assume that the distribution of $Y_{1,1}$ has mean zero and satisfy the following Poincar\'e inequality:
\begin{displaymath}
\forall f \in \CC^1(\R,\R) \ s.t. \ \E(f'(Y_{1,1}))^2 < +\infty , \quad \Var(f(Y_{1,1})) \le \sigma^2 \E(f'(Y_{1,1})^2) \, .
\end{displaymath}
Then, for all deterministic matrices $U \in \MM_n(\R)$ and $V \in \MM_{n,p}(\R)$, and for all integers $n,p$, we have
\begin{equation} \label{concentration1}
\Var \left( \frac1n \Tr(SU) \right) \le \frac{4 \sigma^2 c_n}{n^{5/2}} u(z) \|U\| (\Tr(UU^t))^{1/2}
\end{equation}
and
\begin{multline} \label{concentration2bis}
\Var \left( \frac1n \Tr(X^tSV) \right) \le \frac{9 \sigma^2 c_n}{n^2} v(z) \\
\times \max \left( \|V\| \frac{(\Tr(VV^t))^{1/2}}{n^{1/2}} , \|V\|^{3/2} \frac{(\Tr(VV^t))^{1/4}}{n^{1/4}} , \|V\|^{5/4} \frac{(\Tr(VV^t))^{3/8}}{n^{3/8}} \right) \, ,
\end{multline}
where
\begin{displaymath}
u(z) = \frac{|z|}{|\Imag z|^4} + \frac{1}{|\Imag z|^3}
\end{displaymath}
and
\begin{displaymath}
v(z) = \max \left( \frac{1}{|\Imag z|^2} , \frac{|z|}{|\Imag z|^3} + \frac{1}{|\Imag z|^2} , \left( \frac{|z|}{|\Imag z|^2} + \frac{1}{|\Imag z|} \right)^2 \right) \, .
\end{displaymath}
\end{Proposition}

\begin{Remarks}
\begin{itemize}
\item In the proof of Lemma \ref{freeness_Gaussian_lemma2}, we apply (\ref{concentration1}) to $U=I_n$ and $U=R$, and we apply (\ref{concentration2bis}) to $V=M$.
\item Since $\|V\| \le \Tr(VV^t)^{1/2}$, (\ref{concentration2bis}) implies
\begin{equation} \label{concentration2}
\Var \left( \frac1n \Tr(X^tSV) \right) \le \frac{9 \sigma^2 c_n}{n^{9/4}} v(z) \Tr(VV^t) \, .
\end{equation}
Having in mind the large deviations in Section \ref{section_deviations}, we want to get a bound in $\Tr(MM^t)$ in Lemma \ref{freeness_Gaussian_lemma2}, that is why we use (\ref{concentration2}) instead of (\ref{concentration2bis}) in its proof.
\item We get here slightly better bounds than \cite{VLM}. Indeed, we can recover their results from to ours since $\Tr(AA^t) \le n\|A\|^2$ (see Proposition \ref{prop_trace} (iv)). This improvement is due to the fact that we used the inequality $|\Tr(AC)| \le \sqrt{n} \|A\| (\Tr(CC^*))^{1/2}$ (see Proposition \ref{prop_trace} (iii)) instead of $|\Tr(AC)| \le n \|A\|.\|C\|$.
\item If the distribution of $Y_{1,1}$ satisfies the Poincar\'e inequality with a constant $C$ instead of $\sigma^2$, then $\sigma^2$ must be replaced by $C$ in the bounds (\ref{concentration1}) and (\ref{concentration2bis}).
\item In the case of complex matrices $Y, M, U, V$, the bounds are very similar and only the constants change.
\end{itemize}
\end{Remarks}

\begin{Proof}
Using the sub-additivity property of variance, the Poincar\'e inequality, and the differentiation formula (\ref{DabSjk}), we get
\begin{eqnarray} \label{conc1_eq1}
\lefteqn{   \Var \left( \frac1n \Tr(SU) \right) = \Var \left( \sum_{1 \le j,k \le n} \frac1n S_{j,k}U_{k,j} \right)   } \nonumber \\
& \le & \sigma^2 \E \left[ \sum_{a,b} \left( \frac1n \sum_{j,k} \frac{1}{\sqrt{p}} D_{a,b}S_{j,k}.U_{k,j} \right)^2 \right] \nonumber \\
& = & \sigma^2 \E \left[ \sum_{a,b} \frac{1}{n^2p} \left( \sum_{j,k} (SX)_{j,b}S_{a,k}U_{k,j} + S_{j,a}(X^tS)_{b,k}U_{k,j} \right)^2 \right] \nonumber \\
& = & \frac{\sigma^2}{n^2p} \E \left[ \sum_{a,b} ((SUSX)_{a,b} + (X^tSUS)_{b,a})^2 \right] \nonumber \\
& = & \frac{\sigma^2}{n^2p} \E \left[ \sum_{a,b} (SUSX)_{a,b}^2 + 2(SUSX)_{a,b}(X^tSUS)_{b,a} + (X^tSUS)_{b,a}^2 \right] \nonumber \\
& = & \frac{\sigma^2}{n^2p} \E \left[ \Tr(SUSX(SUSX)^t) + 2 \Tr(SUSXX^tSUS) \right. \nonumber \\
& & \quad \left. + \Tr(X^tSUS(X^tSUS)^t) \right] \, .
\end{eqnarray}

Using the resolvent identity $SXX^t=XX^tS=zS-I_n$, the inequality $|\Tr(AB)| \le \sqrt{n} \Tr(AA^*)^{1/2} \|B\|$ (see Proposition \ref{prop_trace} (iii)), and $\|S\| \le \frac{1}{|\Imag z|}$ (see Proposition \ref{resolvent_covariance} (iv)), we get
\begin{eqnarray} \label{conc1_eq2}
\left| \Tr(SUSX(SUSX)^t) \right| & = & \left| \Tr(U(zS-I_n)SU^tS^2) \right| \nonumber \\
& \le & \sqrt{n} \|U\| \Tr(UU^t)^{1/2} \left( \frac{|z|}{|\Imag z|^4} + \frac{1}{|\Imag z|^3} \right)
\end{eqnarray}
and very similarly,
\begin{equation} \label{conc1_eq3}
\left| \Tr(SUSXX^tSUS) \right| \le \sqrt{n} \|U\| \Tr(UU^t)^{1/2} \left( \frac{|z|}{|\Imag z|^4} + \frac{1}{|\Imag z|^3} \right)
\end{equation}
and
\begin{equation} \label{conc1_eq4}
\left| \Tr(X^tSUS(X^tSUS)^t) \right| \le \sqrt{n} \|U\| \Tr(UU^t)^{1/2} \left( \frac{|z|}{|\Imag z|^4} + \frac{1}{|\Imag z|^3} \right) \, .
\end{equation}

Combining (\ref{conc1_eq1}), (\ref{conc1_eq2}), (\ref{conc1_eq3}) et (\ref{conc1_eq4}), we conclude that
\begin{displaymath}
\Var \left( \frac1n \Tr(SU) \right) \le \frac{4 \sigma^2 c_n}{n^{5/2}} \left( \frac{1}{|\Imag z|^3} + \frac{|z|}{|\Imag z|^4} \right) \|U\| \Tr(UU^t)^{1/2} \, .
\end{displaymath}

Let us now prove the second inequality. By the same arguments as above, we get
\begin{eqnarray} \label{conc2_eq1}
\lefteqn{   \Var \left( \frac1n \Tr(X^tSV) \right) = \Var \left( \sum_{\substack{1 \le k,l \le n \\ 1 \le j \le p}} \frac1n X_{k,j}S_{k,l}V_{l,j} \right)   } \nonumber \\
& \le & \sigma^2 \E \left[ \sum_{a,b} \left( \frac1n \sum_{j,k,l} D_{a,b}X_{k,j} . S_{k,l}V_{l,j} + X_{k,j} . \frac{1}{\sqrt{p}} D_{a,b}S_{k,l} . V_{l,j} \right)^2 \right] \nonumber \\
& = & \frac{\sigma^2}{n^2p} \E \left[ \sum_{a,b} \left( (SV)_{a,b} + (SVX^tSX)_{a,b} + (X^tSVX^tS)_{b,a} \right)^2 \right] \nonumber \\
& = & \frac{\sigma^2}{n^2p} \E \left[ \Tr(SV(SV)^t) + \Tr(SVX^tSX(SVX^tSX)^t) \right. \nonumber \\
& & \quad  + \ \Tr(X^tSVX^tS(X^tSVX^tS)^t) + 2 \Tr(SVX^tSX(SV)^t) \nonumber \\
& & \quad \left. + \ 2 \Tr(SVX^tSVX^tS) + 2 \Tr(SVX^tSXX^tSVX^tS) \right] \, .
\end{eqnarray}

We will now bound these terms, always using the resolvent identities $SXX^t=XX^tS=zS-I_n$, $XX^tS^* = S^*XX^t = \overline{z}S^*-I_n$, inequalities (i)-(iii) in Proposition \ref{prop_trace}, and the bound $\|S\| \le \frac{1}{|\Imag z|}$. We get for example
\begin{equation} \label{conc2_eq2}
\left| \Tr(SV(SV)^t) \right| = \left| \Tr(SVV^tS) \right| \le \sqrt{n} \|V\| \Tr(VV^t)^{1/2} \frac{1}{|\Imag z|^2}
\end{equation}
and
\begin{eqnarray} \label{conc2_eq3}
\lefteqn{   \left| \Tr(SVX^tSX(SVX^tSX)^t) \right|   } \nonumber \\
& = & \left| \Tr(V^tS^2VX^t(zS-I_n)SX) \right| \nonumber \\
& \le & \Tr(V^tS^2VV^t(S^*)^2V)^{1/2} \Tr(X^t(zS-I_n)SXX^tS^*(\overline{z}S^*-I_n)X)^{1/2} \nonumber \\
& = & \Tr(V^tS^2VV^t(S^*)^2V)^{1/2} \Tr(S^*(\overline{z}S^*-I_n)XX^t(zS-I_n)(zS-I_n))^{1/2} \nonumber \\
& \le & \left( \sqrt{n} \|V\|^3 \|S\|^4 \Tr(VV^t)^{1/2} \right)^{1/2} \Tr((\overline{z}S^*-I_n)^2(zS-I_n)^2)^{1/2} \nonumber \\
& \le & n^{1/4} \frac{1}{|\Imag z|^2} \|V\|^{3/2} \Tr(VV^t)^{1/4} . \sqrt{n} \left( \frac{|z|}{|\Imag z|} + 1 \right)^2 \nonumber \\
& = & n^{3/4} \|V\|^{3/2} \Tr(VV^t)^{1/4} \left( \frac{|z|}{|\Imag z|^2} + \frac{1}{|\Imag z|} \right)^2 \, .
\end{eqnarray}
Very similarly, we get
\begin{multline} \label{conc2_eq4}
\left| \Tr(X^tSVX^tS(X^tSVX^tS)^t) \right| \\
\le n^{3/4} \|V\|^{3/2} \Tr(VV^t)^{1/4} \left( \frac{|z|}{|\Imag z|^2} + \frac{1}{|\Imag z|} \right)^2 \, ,
\end{multline}
\begin{equation} \label{conc2_eq5}
\left| \Tr(SVX^tSXV^tS) \right| \le n^{5/8} \|V\|^{5/4} \Tr(VV^t)^{3/8} \left( \frac{|z|}{|\Imag z|^3} + \frac{1}{|\Imag z|^2} \right) \, ,
\end{equation}
\begin{equation} \label{conc2_eq6}
\left| \Tr(SVX^tSVX^tS) \right| \le \sqrt{n} \|V\| \Tr(VV^t)^{1/2} \left( \frac{|z|}{|\Imag z|^3} + \frac{1}{|\Imag z|^2} \right) \, ,
\end{equation}
and finally
\begin{equation} \label{conc2_eq7}
\left| \Tr(SVX^tSXX^tSVX^tS) \right| \le \sqrt{n} \|V\| \Tr(VV^t)^{1/2} \left( \frac{|z|}{|\Imag z|^2} + \frac{1}{|\Imag z|} \right)^2 \, .
\end{equation}

Combining inequalities from (\ref{conc2_eq1}) to (\ref{conc2_eq7}), we finally obtain the announced bound. \qed
\end{Proof}

\begin{Remark}
In the proof above, it is possible to improve some majorizations using the inequality $\|SX\| \le \frac{1}{|\Imag y|}$ where $y$ is a square root of $z$ (see Proposition \ref{resolvent_covariance} (v)). For instance, it allows to get
\begin{displaymath}
|\Tr(X^tSVX^tS(X^tSVX^tS)^t)| \le \sqrt{n} \Tr(VV^t)^{1/2} \|V\| . \frac{1}{|\Imag y|^4}
\end{displaymath}
instead of (\ref{conc2_eq4}). However, in (\ref{conc2_eq3}), which is the other dominant term in (\ref{conc2_eq1}), we can not improve the power of $n$ by this strategy.
\end{Remark}

\subsection{Concentration of the empirical spectral measure}

In Section \ref{section_deviations}, in order to prove Proposition \ref{exponential_equivalence}, we needed the following concentration bound.

\begin{Proposition}[adaptation from {\cite[Theorem 2.5]{BC}}]
Let $\kappa>1$, $Y \in \MM_{n,p}(\R)$ a random matrix with i.i.d. entries bounded by $\kappa$, $M \in \MM_{n,p}(\R)$ a deterministic matrix such that $\frac1n \Tr(MM^t) \le \kappa^2$, and $s,t>0$. We assume that $c_n = \frac{n}{p} \to c \in (0,+\infty)$ as $n \to +\infty$.\\
There exists $\beta>0$ such that for all $s$ large enough, $t$ small enough, $n$ large enough, and $\delta \in \left[ \left( \frac{\beta\kappa^2}{n} \right)^{2/5} , 1 \right]$, we have
\begin{multline} \label{concentration3}
\PR \left( d_{s,t} \left( \mu_{(Y/\sqrt{p}+M)(Y/\sqrt{p}+M)^t} , \E \mu_{(Y/\sqrt{p}+M)(Y/\sqrt{p}+M)^t} \right) \ge \delta \right) \\
\le \frac{\beta\kappa}{\delta^{3/2}} \exp \left( -\frac{n^2\delta^5}{\beta\kappa^4} \right) \, .
\end{multline}
\end{Proposition}

\begin{Remarks}
\begin{itemize}
\item Here $\kappa$ is a constant but we are interested in the dependence on $\kappa$ in the bound since we apply (\ref{concentration3}) to a $\kappa$ depending on $n$ in Section \ref{section_deviations}.
\item This result remains true if $Y$ and $M$ are complex matrices and the entries of $Y$ have independent real and imaginary parts.
\end{itemize}
\end{Remarks}

\begin{Proof}
We will apply \cite[Theorem 1.3(b)]{GZ} to the $(n+p) \times (n+p)$ matrix
\begin{displaymath}
\mathbf{X}_A = \left( \begin{array}{c|c} 0 & \frac{Y}{\sqrt{p}} + M \\ \hline \left( \frac{Y}{\sqrt{p}} + M \right)^t & 0 \end{array} \right) \, .
\end{displaymath}
The matrix $M$ is not present in \cite{GZ} but it is possible to do so because
\begin{displaymath}
\frac{1}{n+p} \Tr(\mathbf{X}_A^2) \le \frac2n \sum_{j,k} \left( \frac{Y_{j,k}}{\sqrt{p}} + M_{j,k} \right)^2 \le \frac{4}{np} \sum_{j,k} Y_{j,k}^2 + \frac4n \Tr(MM^t)
\end{displaymath}
so, thanks to the hypotheses on $Y$ and $M$, we have $\frac{1}{n+p} \Tr(\mathbf{X}_A^2) \le 8 \kappa^2$. Therefore, the argument in \cite[p. 132]{GZ} does not change and we can apply \cite[Theorem 1.3(b)]{GZ} adding the matrix $M$. Consequently, there exists $\beta>0$ such that for all $n$ large enough and $\delta \in \left[ \left( \frac{\beta\kappa^2}{n} \right)^{2/5} , 1 \right]$, we have
\begin{displaymath}
\PR \left( \sup_f \left| \int f \, d\mu_{\mathbf{X}_A} - \E \int f \, d\mu_{\mathbf{X}_A} \right| \ge \delta \right) \le \frac{\beta\kappa}{\delta^{3/2}} e^{-n^2\delta^5/\beta\kappa^4}
\end{displaymath}
where the supremum is taken over all bounded Lipschitz functions $f$ such that
\begin{equation} \label{Lipschitz_norm}
\sup_{x \in \R} |f(x)| + \sup_{x \neq y} \left| \frac{f(x)-f(y)}{x-y} \right| \le 1 \, .
\end{equation}

Moreover, using (\ref{Imy}), we can check that when $s$ is large enough and $t$ is small enough, for every $z \in V_{s,t}$, the function $f : x \mapsto \frac{1}{z-x^2}$ is Lipschitz, bounded, and satisfies (\ref{Lipschitz_norm}). Noting in addition that
\begin{multline*}
\int \frac{1}{z-x^2} \, d\mu_{\mathbf{X}_A}(x) \\
= \frac{1}{n+p} \left( 2n \int \frac{1}{z-x} \, d\mu_{(Y/\sqrt{p}+M)(Y/\sqrt{p}+M)^t}(x) + (n-p).\frac1z \right)
\end{multline*}
and using the definition (\ref{dst}) of $d_{s,t}$, we find (\ref{concentration3}) even if it means to change $\beta$. \qed
\end{Proof}

\section{Technical tools} \label{appendix_tools}

In this appendix, we summarize miscellaneous results used throughout the paper.

\subsection{Traces and matricial norms inequalities} \label{traces_norms}

For a matrix $A \in \MM_{n,p}(\C)$, we denote by $\|A\|$ its operator norm associated to Euclidean norms and
\begin{displaymath}
\|A\|_{\infty} = \max_{1 \le j \le n, 1 \le k \le p} |A_{j,k}| \, .
\end{displaymath}
If $B$ is an other matrix in $\MM_{n,p}(\C)$, we denote by $A \circ B$ the Hadamard product of $A$ and $B$, i.e. the matrix defined by $(A \circ B)_{j,k} = A_{j,k}B_{j,k}$. Finally, $\diag(A)$ denotes the matrix whose entries are given by $A_{j,k}\delta_{j,k}$, where $\delta$ is the Kronecker symbol.

\begin{Proposition} \label{prop_trace}
Let $A,B \in \MM_{n,p}(\C)$, $C \in \MM_{p,n}(\C)$, $D \in \MM_n(\C)$, $E \in \MM_{p,q}(\C)$. We have the following.
\begin{itemize}
\item[(i)] $|\Tr(AC)| \le (\Tr(AA^*))^{1/2} (\Tr(CC^*))^{1/2}$,
\item[(ii)] $|\Tr(AC)| \le n \|A\|.\|C\|$,
\item[(iii)] $|\Tr(AC)| \le \sqrt{n} \|A\| (\Tr(CC^*))^{1/2}$,
\item[(iv)] $\|A\| \le (\Tr(AA^*))^{1/2} \le \sqrt{n} \|A\|$,
\item[(v)] $(\Tr(AA^*))^{1/2} \le \sqrt{np} \|A\|_{\infty}$,
\item[(vi)] $\| \diag(D) \| = \| \diag(D) \|_{\infty} \le \|D\|_{\infty}$,
\item[(vii)] $\| A \circ B \|_{\infty} \le \| A \|_{\infty} \| B \|_{\infty}$,
\item[(viii)] $\| A \circ B \| \le \|A\|.\|B\|$,
\item[(ix)] $\Tr((A \circ B)(A \circ B)^*) \le \Tr(AA^*) \|B\|_{\infty}^2$.
\end{itemize}
\end{Proposition}

Most of these points are classical or easy to check. Note that the combination (iii) of (i) and (ii) will be crucial for us and that a proof of (viii) requires the use of the matrices
\begin{displaymath}
A' = \left( \begin{array}{ccccccccccc} A_{1,1} & & & A_{1,2} & & & & A_{1,p} & & \\ & \ddots & & & \ddots & & \ldots & & \ddots & \\ & & A_{n,1} & & & A_{n,2} & & & & A_{n,p} \end{array} \right)
\end{displaymath}
and
\begin{displaymath}
B' = \left( \begin{array}{cccc} B_{1,1} \\ \vdots \\ B_{n,1} \\ & B_{1,2} \\ & \vdots \\ & B_{n,2} \\ & & \ddots \\ & & & B_{1,p} \\ & & & \vdots \\ & & & B_{n,p} \end{array} \right) \, .
\end{displaymath}

\subsection{Properties of resolvents} \label{resolvents}

Let $A \in \HH_n(\C)$ and $z \in \C \setminus \R$. The \emph{resolvent of $A$ at $z$} is the matrix $R(A) = (zI_n-A)^{-1}$. For $A \in \MM_{n,p}(\C)$, we denote by $S(A)$ the resolvent $R(AA^*)$, or just $S$ if no confusion can arise.

\begin{Proposition} \label{resolvent_covariance}
Let $A,B \in \MM_{n,p}(\C)$ and $z \in \C \setminus \R$. We have the following.
\begin{itemize}
\item[(i)] $SAA^* = AA^*S = zS - I_n$,
\item[(ii)] $S(A+B)-S(A) = S(A+B)(AB^*+BA^*+BB^*)S(A)$,
\item[(iii)] $G_{\mu_{AA^*}}(z) = \frac1n \Tr S$,
\item[(iv)] $\|S\|_{\infty} \le \|S\| \le \frac{1}{|\Imag z|}$,
\item[(v)] $\|SA\|_{\infty} \le \|SA\| \le \frac{1}{|\Imag y|}$, where $y$ is a square root of $z$,
\item[(vi)] $\|A^*SA\|_{\infty} \le \|A^*SA\| \le 1 + \left| \frac{z}{\Imag z} \right|$.
\item[(vii)] We denote by $D_{a,b}$ the derivation w.r.t. $\Reel A_{a,b}$ and $\delta$ the Kronecker symbol. For all $a,j,k \in \llbracket 1,n \rrbracket$ and $b,l,m \in \llbracket 1,p \rrbracket$, we have
\begin{equation} \label{DabSjk}
D_{a,b}S_{j,k} = (SA)_{j,b}S_{a,k} + S_{j,a}(A^*S)_{b,k} \, ,
\end{equation}
\begin{equation} \label{Dab(SX)jk}
D_{a,b}(SA)_{j,l} = (SA)_{j,b}(SA)_{a,l} + S_{j,a}(A^*SA)_{b,l} + \delta_{b,l} S_{j,a} \, ,
\end{equation}
\begin{equation} \label{Dab(X^*S)jk}
D_{a,b}(A^*S)_{l,k} = (A^*SA)_{l,b}S_{a,k} + (A^*S)_{l,a}(A^*S)_{b,k} + \delta_{b,l} S_{a,k} \, ,
\end{equation}
\begin{multline} \label{Dab(X^*SX)jk}
D_{a,b}(A^*SA)_{l,m} = (A^*SA)_{l,b}(SA)_{a,m} + (A^*S)_{l,a}(A^*SA)_{b,m} \\
+ \, \delta_{b,m} (A^*S)_{l,a} + \delta_{b,l} (SA)_{a,m} \, ,
\end{multline}
\begin{multline} \label{D2abSjk}
D^2_{a,b}S_{j,k} = 2 [S_{j,a}S_{a,k} + (SA)_{j,b}(SA)_{a,b}S_{a,k} + S_{j,a}(A^*SA)_{b,b}S_{a,k} \\
+ \, S_{j,a}(A^*S)_{b,a}(A^*S)_{b,k} + (SA)_{j,b}S_{a,a}(A^*S)_{b,k}] \, ,
\end{multline}
\begin{multline} \label{D3abSjk}
D^3_{a,b}S_{j,k} = 6 [(SA)_{j,b}S_{a,a}S_{a,k} + S_{j,a}(A^*S)_{b,a}S_{a,k} + S_{j,a}(SA)_{a,b}S_{a,k} \\
+ S_{j,a}S_{a,a}(A^*S)_{b,k} + (SA)_{j,b}(SA)_{a,b}^2S_{a,k} + S_{j,a}(A^*S)_{b,a}^2(A^*S)_{b,k} \\
+ (SA)_{j,b}(SA)_{a,b}S_{a,a}(A^*S)_{b,k} + (SA)_{j,b}S_{a,a}(A^*S)_{b,a}(A^*S)_{b,k} \\
+ S_{j,a}(A^*SA)_{b,b}(SA)_{a,b}S_{a,k} + (SA)_{j,b}S_{a,a}(A^*SA)_{b,b}S_{a,k} \\
+ S_{j,a}(A^*S)_{b,a}(A^*SA)_{b,b}S_{a,k} + S_{j,a}(A^*SA)_{b,b}S_{a,a}(A^*S)_{b,k}] \, .
\end{multline}
\end{itemize}
\end{Proposition}

Most of these relations are classical or obtained by simple computations. Note however that (v) and (vi) respectively follow from the identities
\begin{displaymath}
\left( \begin{array}{c|c} \overline{y}I_n & -A \\ \hline -A^* & \overline{y}I_p \end{array} \right)^{-1} = \left( \begin{array}{c|c} \overline{y} S^* & A(\overline{y}^2 I_p-A^*A)^{-1} \\ \hline (SA)^* & \overline{y} (\overline{y}^2 I_p-A^*A)^{-1} \end{array} \right)
\end{displaymath}
and $A^*(zI_n-AA^*)^{-1}A = -I_p + z(zI_p-A^*A)^{-1}$.

Note also that if $y$ is a square root of $z$ and $z$ belongs to the domain $V_{s,t}$ defined by (\ref{Vst}), then we can easily prove that
\begin{equation} \label{Imy}
|\Imag y| = \frac{\Imag z}{2} \left( \sqrt{\frac{(\Reel z)^2}{(\Imag z)^2} + 1} - \frac{\Reel z}{\Imag z} \right) > \left( \frac{s}{2} (\sqrt{t^2+1}-t) \right)^{1/2} > 0 \, .
\end{equation}

\subsection{Inequalities for empirical spectral measures} \label{inequalities_ESM}

\begin{Proposition}[Rank inequality, see {\cite[Lemma B.1]{BC}}]
Let $A,B \in \HH_n(\C)$. We have
\begin{equation} \label{rank_inequality}
\dKS(\mu_A,\mu_B) \le \frac1n \rg(A-B) \, .
\end{equation}
\end{Proposition}

\begin{Proposition}[Rank inequality for covariance matrices, see {\cite[Theorem A.44]{BS}}]
Let $A,B \in \MM_{n,p}(\C)$. We have
\begin{equation} \label{rank_covariance}
\dKS(\mu_{AA^*},\mu_{BB^*}) \le \frac1n \rg(A-B) \, .
\end{equation}
\end{Proposition}

\begin{Proposition}[Hoffman-Wielandt inequality, see {\cite[Lemma B.2]{BC}}]
Let $A,B \in \HH_n(\C)$. We have
\begin{equation} \label{HW_inequality}
W_2^2(\mu_A,\mu_B) \le \frac1n \Tr((A-B)^2) \, ,
\end{equation}
where $W_2$ denotes the $L^2$-Wasserstein distance on $\PP(\R)$.
\end{Proposition}

\begin{Proposition}[see {\cite[Corollary A.42]{BS}}]
Let $A,B \in \MM_{n,p}(\C)$. We have
\begin{equation} \label{HW_covariance}
W_2^4(\mu_{AA^*},\mu_{BB^*}) \le \frac{2}{n^2} \Tr(AA^*+BB^*) \Tr((A-B)(A-B)^*) \, .
\end{equation}
\end{Proposition}

\begin{Proposition}[Schatten's inequality, see {\cite[Theorem 3.32]{Zhan}}]
Let $A \in \HH_n(\C)$ and $p \in (0,2]$. We have
\begin{equation} \label{Schatten}
\int_{\R} |x|^p \, d\mu_A(x) \le \frac1n \sum_{k=1}^n \left( \sum_{j=1}^n |A_{k,j}|^2 \right)^{p/2} \, .
\end{equation}
\end{Proposition}

\section*{Acknowledgements}

I gratefully thank my supervisors Catherine Donati-Martin and Myl\`ene Ma\"ida, not only for the fruitful discussions we shared around this work, but also for their constant support and their availability. I also thank Charles Bordenave for the discussions we had especially during the conference "\'Etats de la Recherche en Matrices Al\'eatoires" organized in December 2014 in IHP, Paris.

\bibliographystyle{plain}
\bibliography{article1_biblio}

\end{document}